\documentclass{article}
\usepackage{graphicx} 
\usepackage{amsmath}
\usepackage{amsthm}
\usepackage{amssymb}
\usepackage{amscd}
\usepackage{amsbsy}
\usepackage{amsfonts}
\usepackage{euscript}
\usepackage{stmaryrd}
\usepackage{mathrsfs}
\usepackage[matrix,arrow,curve,all,2cell,cmtip]{xy}

\newtheorem{theorem}{Theorem}
\newtheorem{lemma}{Lemma}
\newtheorem{proposition}{Proposition}
\newtheorem{corollary}{Corollary}
\newtheorem{definition}{Definition}
\newtheorem{remark}[theorem]{Remark}
\newtheorem{example}{Example}

\title{Double factorization systems in equivariant topology and topos theory}

\author{E.\,V.~Martyanov}

\date{}

\begin{document}

\maketitle

\begin{abstract}
In the present work, we investigate the extension of double factorization systems to the categories of Eilenberg-Moore (co)algebras. We show that the double factorization systems $(\texttt{ExEpi},\texttt{Bim},\texttt{ExMono})$ in the categories $\mathbf{Tych}$, $\mathbf{Unif}$ and $\mathbf{Comp}$ extend to the same double factorization systems in the corresponding categories of the Eilenberg-Moore algebras $\mathbf{Tych}^{\mathbb{H}^t}$, $\mathbf{Unif}^{\mathbb{H}^u}$ and $\mathbf{Comp}^{\mathbb{H}^c}$. We establish a connection between cartesian double factorization systems and LT-topologies. We provide sufficient conditions for the extension of cartesian double factorization systems to the topos of coalgebras.

\


\end{abstract}

\section{Introduction}
Since the inception of category theory, one of its key tasks has been the factorization of morphisms, i.e. their representation as a composition of morphisms from certain classes. For a continuous mapping $f\colon X\to Y$ from a topological space $X$ to a topological space $Y$ there exists a representation as a composition $hg$, where $g$ is a continuous mapping onto the image $f(X)$, and $h\colon f(X)\hookrightarrow Y$ is the embedding of the image $f(X)$ into the space $Y$. In turn, the mapping $g$ can be represented as a composition of a continuous quotient mapping and a continuous bijection. Thus, every continuous mapping can be expressed as a composition of a continuous quotient map, a continuous bijection and an embedding. A similar representation holds for uniformly continuous maps of uniform spaces. In this case, an arbitrary uniformly continuous map $f\colon(X,\mathcal{U})\to(Y,\mathcal{V})$ from a uniform space $(X,\mathcal{U})$ to a uniform space $(Y,\mathcal{V})$ can be represented as a composition of a uniform quotient map, a dense uniformly continuous injection and a closed uniformly continuous embedding. In the terminology of category theory, this means that the categories of topological spaces $\mathbf{Top}$ and uniform spaces $\mathbf{Unif}$ have a double factorization system (see Section 4.1).

In any Grothendieck topos, there exists a factorization system $(\texttt{Epi},\texttt{Mono})$. If a Grothendieck topology is defined on a Grothendieck topos, then every monomorphism admits further factorization into a closed and a dense monomorphism. In other words, the subcategory of monomorphisms $\texttt{Mono}$ in a Grothendieck topos has a factorization system $(\texttt{DnsMono},\texttt{ClsMono})$. A similar factorization holds in an arbitrary topos if a Lawvere-Tierney topology (LT-topology) is defined on it.

In the examples above, the morphisms of the categories admit a representation as a composition of morphisms from three classes $\mathcal{E}$, $\mathcal{J}$ and $\mathcal{M}$. In \cite{PT}, the concept of a double factorization system $(\mathcal{E},\mathcal{J},\mathcal{M})$ was introduced in a category $\mathcal{C}$. A double factorization system (dfs) is a direct generalization of the notion of a factorization system (Definition \ref{d3}). The presence of such a dfs is an important property of a category. For instance, certain dfs's on a topos allow us to define LT-topologies and, consequently, to define reflective subcategories of sheaves.

The problem of extending such structures is studied in the theory of model categories (basic concepts of the model category theory can be found in \cite{H}). Suppose that an adjunct $F\colon\mathcal{C}\leftrightarrows\mathcal{D}\colon G$ is given from a bicomplete category $\mathcal{C}$ to a model category $(\mathcal{D},\mathcal{C}of_{\mathcal{D}},\mathcal{W}_{\mathcal{D}},\mathcal{F}ib_{\mathcal{D}})$. Is it possible to transfer the model structure $(\mathcal{C}of_{\mathcal{D}},\mathcal{W}_{\mathcal{D}},\mathcal{F}ib_{\mathcal{D}})$ from the category $\mathcal{D}$ to the category $\mathcal{C}$ in such a way that the given adjunction becomes a Quillen adjunction? The most natural way to perform such a transfer is to define the class of cofibrations and the class of weak equivalences in $\mathcal{C}$ as the preimages of the corresponding classes in $\mathcal{D}$, i.e. $\mathcal{C}of_{\mathcal{C}}=F^{-1}(\mathcal{C}of_{\mathcal{D}})$ and $\mathcal{W}_{\mathcal{C}}=F^{-1}(\mathcal{W}_{\mathcal{D}})$. The class of fibrations is defined as the class of all morphisms of the category $\mathcal{C}$ that have the right lifting property (Definition \ref{d1}) with respect to trivial cofibrations, i.e. $\mathcal{F}ib_{\mathcal{C}}=(F^{-1}(\mathcal{C}of_{\mathcal{D}})\cap F^{-1}(\mathcal{W}_{\mathcal{D}}))^\boxslash$ (see \cite{BH}). The resulting system may not satisfy the axioms of a model category. One of the obstacles is the impossibility of representing a morphism as a composition of a fibration and a trivial cofibration (resp., a trivial fibration and a cofibration). To solve this problem, the small object argument is used, which in turn relies on the theory of ordinals. However, imposing an additional condition of uniqueness in the definition of the right (left) lifting property allows us to avoid using the theory of ordinals to solve the problem of extending a dfs.

The requirement of uniqueness in the definition of the left (resp., right) lifting property leads to the notion of an (orthogonal) Quillen factorization system (Definition 3.9 \cite{PT}). If a Quillen factorization system $(\mathcal{C}of_{\mathcal{C}},\mathcal{W}_{\mathcal{C}},\mathcal{F}ib_{\mathcal{C}})$ is given on a (finitely) bicomplete category $\mathcal{C}$, then the model category $(\mathcal{C},\mathcal{C}of_{\mathcal{C}},\mathcal{W}_{\mathcal{C}},\mathcal{F}ib_{\mathcal{C}})$ will be called \textit{orthogonal}. There exists a connection between dfs's and Quillen factorization systems. Every Quillen factorization system on a (finitely) bicomplete category defines a dfs (Theorem 3.10 \cite{PT}). The converse statement is not true in general (see Example 3.12 \cite{PT}).

The present paper is devoted to the extension of dfs's in the categories $\mathbf{Tych}$, $\mathbf{Unif}$ and $\mathbf{Comp}$ to the categories of Eilenberg-Moore algebras $\mathbf{Tych}^{\mathbb{H}^t}$, $\mathbf{Unif}^{\mathbb{H}^u}$ and $\mathbf{Comp}^{\mathbb{H}^c}$, respectively, as well as to the extension of dfs's to the topos of coalgebras. We will proof the possibility of extending dfs's of the form $(\texttt{ExEpi},\texttt{Bim},\texttt{ExMono})$ in the categories $\mathbf{Unif}$, $\mathbf{Tych}$ and $\mathbf{Comp}$ to analogous dfs's in the categories $\mathbf{Tych}^{\mathbb{H}^t}$, $\mathbf{Unif}^{\mathbb{H}^u}$ and $\mathbf{Comp}^{\mathbb{H}^c}$. For toposes, we will establish a connection between cartesian dfs's (Definition \ref{d16}) and LT-topologies.

In \textbf{Section 2}, brief information and necessary definitions related to dfs and model categories are provided. The main results of this section are Proposition \ref{p1} and Proposition \ref{p2}. Let $\mathbb{T}=(T,\eta,\mu)$ be a monad on a (finitely) bicomplete category $\mathcal{C}$. There exists an adjunction $F^{\mathbb{T}}\colon\mathcal{C}\leftrightarrows\mathcal{C}^{\mathbb{T}}\colon V^{\mathbb{T}}$, in which the right adjoint $V^{\mathbb{T}}$ is the forgetful functor, and the left adjoint $F^{\mathbb{T}}$ is the free algebra functor, assigning to each object of the category $\mathcal{C}$ an Eilenberg-Moore algebra. If the underlying functor $T$ of the monad $\mathbb{T}$ satisfies the conditions $T(\mathcal{E})\subseteq\mathcal{E}$ (preserves trivial cofibrations) and $T(\mathcal{J})\subseteq\mathcal{J}$ (preserves bifibrant morphisms), then the forgetful functor $V^{\mathbb{T}}$ generates the right-induced dfs (Definition \ref{d12}) in the category of Eilenberg-Moore algebras $\mathcal{C}^{\mathbb{T}}$. Moreover, with respect to the right-induced dfs in the category of Eilenberg-Moore algebras and the dfs on $\mathcal{C}$, the forgetful functor $V^{\mathbb{T}}$ is not only a right Quillen functor (Definition \ref{d9}), but also preserves and reflects the corresponding local objects (Proposition \ref{p1}). By duality principle, the results formulated in Proposition \ref{p1} can be transferred to the category of Eilenberg-Moore coalgebras.

In topos theory, there is a connection between the topos of coalgebras and the category of sheaves (see \cite{MM}). Since sheaves are $\mathcal{J}$-local objects for a suitable dfs $(\mathcal{E},\mathcal{J},\mathcal{M})$, it is natural to consider the relationship between the category of Eilenberg-Moore coalgebras $\mathcal{C}_{\mathbb{G}}$ and the category of $\mathcal{J}$-local objects generated by the dfs $(\mathcal{E},\mathcal{J},\mathcal{M})$ on the category $\mathcal{C}$ (Proposition \ref{p2}).

In \textbf{Section 3}, we investigate the connection between dfs's and LT-topologies on an arbitrary topos. The main results of this section are Theorem \ref{th7} and Proposition \ref{p4}. Every LT-topology corresponds to a dfs $(\texttt{Epi}, \texttt{DnsMono}, \texttt{ClsMono})$. This dfs is stable under pullbacks. Conversely, every cartesian dfs generates a certain LT-topology (Proposition \ref{p3}).

Theorem \ref{th7} establishes the connection between Quillen adjunctions and cartesian dfs's. Let $(\mathcal{E}_1,\mathcal{J}_1,\mathcal{M}_1)$ and $(\mathcal{E}_2,\mathcal{J}_2,\mathcal{M}_2)$ be cartesian dfs's on toposes $\mathcal{C}_1$ and $\mathcal{C}_2$, respectively. By Proposition \ref{p3} $(\mathcal{E}_1,\mathcal{J}_1,\mathcal{M}_1)$ and $(\mathcal{E}_2,\mathcal{J}_2,\mathcal{M}_2)$ generates LT-topologies $k_1$ and $k_2$, respectively. If $F\colon\mathcal{C}_1\leftrightarrows\mathcal{C}_2\colon G$ is a Quillen adjunction and $F$ preserves finite limits, then
$G$ maps $k_2$-sheaves (resp. $k_2$-separating objects) to $k_1$-sheaves (resp. $k_1$-separating objects). Moreover, item (b) of Theorem \ref{th7} establishes equivalent conditions for preservation of cofibrations by the functor $G$.

If the underlying functor $G$ of a cartesian comonad $\mathbb{G} = (G,\epsilon,\delta)$ on the topos $\mathcal{C}$ preserves bifibrant morphisms and trivial fibrations, and is $(k,k)$-continuous (Definition in 5.7. \cite{DT} and Remark \ref{r8}), then the LT-topology $k_{\mathbb{G}}$ is an extension $\widetilde{k}$ of the LT-topology $k$ (Proposition \ref{p4}).

In \textbf{Section 4}, we investigate the extension of dfs to categories from equivariant topology $\mathbf{Tych}^{\mathbb{H}^t}$, $\mathbf{Unif}^{\mathbb{H}^u}$ and $\mathbf{Comp}^{\mathbb{H}^c}$ (see Section 4.2). Even for the case of factorization systems in these categories, obstacles arise for the existence of their extensions. Let $me$ be the factorization of an equivariant mapping $f\colon(X,\alpha)\to(Y,\beta)$. The extension of mappings $m$ and $e$ to equivariant ones is possible if and only if there exists a corresponding action $\beta'\colon G\times Y'\to Y'$ satisfying the equalities $\beta(1_G\times m)=m\beta'$ and $\beta'(1_G\times e)=e\alpha$ (see Section 4.3.1 \cite{V}).

The main results of this chapter are Proposition \ref{p5} and Theorem \ref{th8}. One of the most important properties of a category is cocompleteness. Proposition \ref{p5} proves the cocompleteness of $\mathbf{Unif}^{\mathbb{H}^u}$, $\mathbf{Tych}^{\mathbb{H}^t}$ and $\mathbf{Comp}^{\mathbb{H}^c}$. Theorem \ref{th8} formulates sufficient conditions for extension of the dfs $(\texttt{ExEpi},\texttt{Bim},\texttt{ExMono})$ to the categories $\mathbf{Tych}^{\mathbb{H}^t}$, $\mathbf{Unif}^{\mathbb{H}^u}$ and $\mathbf{Comp}^{\mathbb{H}^c}$. It is also shown that every object in $\mathbf{Tych}^{\mathbb{H}^t}$, $\mathbf{Unif}^{\mathbb{H}^u}$ and $\mathbf{Comp}^{\mathbb{H}^c}$ is $\texttt{Epi}$-local if and only if it is $\texttt{Epi}$-local in $\mathbf{Tych}$, $\mathbf{Unif}$ and $\mathbf{Comp}$, respectively.

The categories $\mathbf{Tych}^{\mathbb{H}^t}$, $\mathbf{Unif}^{\mathbb{H}^u}$ and $\mathbf{Comp}^{\mathbb{H}^c}$ are connected by adjunctions $F^{tu}\colon\mathbf{Tych}^{\mathbb{H}^t}\leftrightarrows\mathbf{Unif}^{\mathbb{H}^u}\colon T^{tu}$, $\beta^{tc}\colon\mathbf{Tych}^{\mathbb{H}^t}\leftrightarrows\mathbf{Comp}^{\mathbb{H}^c}\colon\iota^{tc}$ and $L^{uc}\colon\mathbf{Unif}^{\mathbb{H}^u}\leftrightarrows\mathbf{Comp}^{\mathbb{H}^c}\colon R^{uc}$ (Theorem \ref{th8}). Moreover, these adjunctions are obtained as corresponding extensions of adjunctions $F\colon\mathbf{Tych}\leftrightarrows\mathbf{Unif}\colon T$, $\beta\colon\mathbf{Tych}\leftrightarrows\mathbf{Comp}\colon\iota$ and $L^u\colon\mathbf{Unif}\leftrightarrows\mathbf{Comp}\colon R^u$, where $R^u=i^uF^c$, $L^u=T^c\beta^u$, $\beta$ is the Stone-\v{C}ech compactification functor, $\beta^u$ is the Samuel compactification functor, $\iota$ and $\iota^u$ are the corresponding embedding functors, $T^c$ and $F^c$ are functors obtained by restricting to subcategories $\mathbf{Comp}$ and $\mathbf{CBUnif}$ the topologization functor $T$ and the finest uniformity functor $F$.

\section{Categorical preliminaries}

\subsection{Double factorization systems and Quillen factorization systems.}

We assume that the reader is familiar with basic definitions of category theory. All necessary information can be found in standard monographs (see, for example, \cite{B} or \cite{M}). All categories in the text will be considered \textit{large}, meaning that objects and morphisms of these categories are subsets (=``classes``) of some universe $\mathfrak{U}$. By \textit{small} categories we mean categories that are elements (=``sets``) of this universe $\mathfrak{U}$. Unless otherwise stated, we will assume that the category $\mathcal{C}$ is bicomplete.

\begin{definition}\label{d1}
\rm{Let $\mathcal{C}$ be a category. A morphism $f$ in $\mathcal{C}$ is said to have the \textit{left lifting property} with respect to a morphism $g$ if for any diagram $$
\xymatrix{
  .\ar[d]_f \ar[r] & .\ar[d]^{g} \\
  . \ar@{-->}[ru]^{h} \ar[r] & .   \\
}
$$
there exists a solution $h$ making both triangles commute. It is also said that $g$ has the \textit{right lifting property} with respect to $f$. (\emph{Notation: $f\boxslash g$}).

If the solution $h$ is unique, then $f$ is said to be \textit{left orthogonal} to $g$ (respectively, $g$ is \textit{right orthogonal} to $f$), and the notation $f\bot g$ is used.}
\end{definition}

Let $M$ be an arbitrary class of morphisms in category $\mathcal{C}$. Then we define: $$M^{\boxslash}=\{f\mid \forall g(g\in M\Rightarrow g\boxslash f)\}$$
$$^{\boxslash}M=\{f\mid \forall g(g\in M\Rightarrow f\boxslash g)\}.$$ The notations $M^{\bot}$ and $^{\bot}M$ have analogous meanings.

Let $\mathcal{R}$ and $\mathcal{L}$ be two classes of morphisms in category $\mathcal{C}$. We denote by $\mathcal{R} \cdot \mathcal{L}$ the class of morphisms consisting of all possible compositions of the form $me$, where $m \in \mathcal{R}$ and $e \in \mathcal{L}$. We call the classes $\mathcal{L}$ and $\mathcal{R}$ \textit{mutually orthogonal} (notation $\mathcal{L}\perp\mathcal{R}$) if for any $e\in\mathcal{L}$ and $m\in\mathcal{R}$ the condition $e\perp m$ holds.

\begin{definition}\label{d2}
\rm{A \emph{weak factorization system} (\emph{wfs}) on a category $\mathcal{C}$ is an ordered pair $(\mathcal{L}, \mathcal{R})$ consisting of two classes of morphisms in $\mathcal{C}$ satisfying the following conditions:}
\begin{itemize}
\item [(i)] $\texttt{Mor}\ \mathcal{C} = \mathcal{R} \cdot \mathcal{L}$;
\item [(ii)] $\mathcal{L} = {}^{\boxslash}\mathcal{R}$ and $\mathcal{R} = \mathcal{L}^{\boxslash}$.
\end{itemize}
\end{definition}

\begin{definition}\label{d3}
\rm{A \emph{factorization system} (\emph{fs}) is a weak factorization system $(\mathcal{L}, \mathcal{R})$ such that
$\mathcal{L} = {}^{\bot}\mathcal{R}$ and $\mathcal{R} = \mathcal{L}^{\bot}$}.
\end{definition}

Let $(\mathcal{L}, \mathcal{R})$ be an fs in category $\mathcal{C}$. Then the following conditions hold (cf. Definition 5.5.1 Vol. 1 \cite{B}):
\begin{itemize}
  \item [(i)] $\texttt{Iso}\subseteq\mathcal{L}\cap\mathcal{R}$;
  \item [(ii)] $\mathcal{L}\cdot\mathcal{L}\subseteq\mathcal{L}$ and $\mathcal{R}\cdot\mathcal{R}\subseteq\mathcal{R}$;
  \item [(iii)] $\mathcal{L}\perp\mathcal{R}$.
\end{itemize}

In \cite{PT} was introduced the concept of a double factorization system.

\begin{definition}[Definition 2.1 \cite{PT}]\label{d4}
\rm{A \emph{double} (\emph{orthogonal}) \emph{factorization system} (\emph{dfs}) in a category $\mathcal{C}$ is an ordered triple
$(\mathcal{E}, \mathcal{J}, \mathcal{M})$ of morphism classes satisfying the following conditions}:
\begin{enumerate}
\item [(i)] $\texttt{Iso}\ \mathcal{C}\cdot\mathcal{E}\subseteq\mathcal{E}$ and $\texttt{Iso}\ \mathcal{C}\cdot\mathcal{J}\cdot\texttt{Iso}\ \mathcal{C}\subseteq\mathcal{J}$ and $\mathcal{M}\cdot\texttt{Iso}\ \mathcal{C}\subseteq\mathcal{M}$;
\item [(ii)] $\texttt{Mor}\ \mathcal{C}=\mathcal{M}\cdot\mathcal{J}\cdot\mathcal{E}$;
\item [(iii)] For any commutative diagram in $\mathcal{C}$ of the form:
$$
\xymatrix{
  A \ar[d]_{u} \ar[r]^{e} & B \ar@{.>}[dl]_s \ar[r]^{j} & C \ar@{.>}[dl]^t \ar[d]^{v} \\
  D \ar[r]_{j'} & E \ar[r]_{m} & F   }
$$
where $e \in \mathcal{E}$, $j, j' \in \mathcal{J}$, $m \in \mathcal{M}$, there exists a unique pair of arrows $(s, t)$ such that $se = u$, $j's = tj$ and $mt = v$.
\end{enumerate}
\end{definition}

Let $(\mathcal{E}, \mathcal{J}, \mathcal{M})$ be a dfs in category $\mathcal{C}$. Then the following conditions hold (see \cite{PT} for details):
\begin{enumerate}
    \item[(i)] $\mathcal{E} ={}^\perp(\mathcal{M} \cdot \mathcal{J})$, $\mathcal{J} ={}^\perp\mathcal{M}\cap\mathcal{E}^\perp$, $\mathcal{M} = (\mathcal{J} \cdot \mathcal{E})^\perp$;

    \item[(ii)] $(\mathcal{E}, \mathcal{J}, \mathcal{M})$-factorization of morphisms is unique up to isomorphism;

    \item[(iii)] $\texttt{Iso}\ \mathcal{C} \subseteq \mathcal{E} \cap \mathcal{J} \cap \mathcal{M}$.
\end{enumerate}

\begin{definition}\label{d5}
\rm{A (\emph{orthogonal}) \emph{Quillen factorization system} (\emph{qfs}) on a category $\mathcal{C}$ is given by classes of morphisms $\mathcal{F}ib$, $\mathcal{W}$, $\mathcal{C}of$ satisfying}:
\begin{itemize}
\item[(i)] $(\mathcal{W} \cap \mathcal{C}of, \mathcal{F}ib)$ and $(\mathcal{C}of, \mathcal{W} \cap \mathcal{F}ib)$ are factorization systems;
\item[(ii)] $\mathcal{W}$ has the \emph{2-out-of-3 property}, i.e. for any two composable morphisms $f$ and $g$ in $\mathcal{W}$, if two of $f$, $g$ and $gf$ are in $\mathcal{W}$, then so is the third.
\end{itemize}
\end{definition}

A morphism $f$ in a category $\mathcal{C}$ with a qfs $(\mathcal{C}of,\mathcal{W},\mathcal{F}ib)$ is called:
\begin{itemize}
\item a \textit{cofibration} if $f\in\mathcal{C}of$;
\item a \textit{fibration} if $f\in\mathcal{F}ib$;
\item a \textit{bifibrant morphism} if $f\in\mathcal{B}if=\mathcal{F}ib\cap\mathcal{C}of$;
\item a \textit{weak equivalence} if $f\in\mathcal{W}$;
\item a \textit{trivial fibration} if $f\in\mathcal{F}ib\cap\mathcal{W}$;
\item a \textit{trivial cofibration} if $f\in\mathcal{C}of\cap\mathcal{W}$.
\end{itemize}

An object $X$ in a category $\mathcal{C}$ with an qfs $(\mathcal{C}of,\mathcal{W},\mathcal{F}ib)$ is called:
\begin{itemize}
\item \textit{fibrant} if the morphism $X\to 1$ is a fibration, where 1 is the terminal object of $\mathcal{C}$;
\item \textit{cofibrant} if the morphism $0\to X$ is a cofibration, where 0 is the initial object of $\mathcal{C}$;
\item \textit{bifibrant} if $X$ is both fibrant and cofibrant.
\end{itemize}

We denote by $\mathcal{C}_c$ (respectively $\mathcal{C}_f$, $\mathcal{C}_{cf}$) \textit{the full subcategory of cofibrant} (respectively \textit{fibrant}, \textit{bifibrant}) objects of $\mathcal{C}$. The connection between dfs and qfs is given in the following statement.

\begin{theorem}[Theorem 3.10 \cite{PT}]\label{th1}
For any dfs $(\mathcal{E}, \mathcal{J}, \mathcal{M})$ satisfying:
\begin{itemize}
\item[(i)] $\mathcal{E} \cdot \mathcal{M} \subseteq \mathcal{M} \cdot \mathcal{E}$,
\item[(ii)] $j \in \mathcal{J}$ is an isomorphism if there exists $e \in \mathcal{E}$ such that $ej \in \mathcal{E}$, or $m \in \mathcal{M}$ such that $jm \in \mathcal{M}$,
\end{itemize}
the triple $(\mathcal{J} \cdot \mathcal{E}, \mathcal{M} \cdot \mathcal{E}, \mathcal{M} \cdot \mathcal{J})$ forms an qfs. Conversely, for any qfs $(\mathcal{C}of, \mathcal{W},\mathcal{F}ib)$, the triple $(\mathcal{C}of\cap\mathcal{W}, \mathcal{B}if, \mathcal{F}ib\cap\mathcal{W})$ forms a dfs $(\mathcal{E}, \mathcal{J}, \mathcal{M})$ satisfying conditions \emph{(i)} and \emph{(ii)}. These correspondences establish a bijection between all qfss of the category $\mathcal{C}$ and all dfss of $\mathcal{C}$ satisfying conditions \emph{(i)} and \emph{(ii)}.
\end{theorem}

\begin{remark}\label{r2}
{\rm There exist dfs that do not correspond to any qfs. For example, in the dfs $(\texttt{RegEpi}, \texttt{Bim}, \texttt{RegMono})$, the class of weak equivalences $\mathcal{W} = \texttt{RegMono} \cdot \texttt{RegEpi}$ rarely satisfies the 2-out-of-3 property, except when $\texttt{Bim}=\texttt{Iso}$ (see details in Example 3.12 (2) \cite{PT}).}
\end{remark}

\begin{remark}[Definition 3.1 \cite{PT}]\label{r3}
{\rm For ofs we use the same terminology as for qfs:
\begin{itemize}
\item $\mathcal{E}=\text{trivial cofibrations}$;
\item $\mathcal{J}=\text{bifibrant morphisms}$;
\item $\mathcal{M}=\text{trivial fibrations}$;
\item $\mathcal{J}\cdot\mathcal{E}=\text{cofibrations}$;
\item $\mathcal{M}\cdot\mathcal{E}=\text{weak equivalences}$;
\item $\mathcal{M} \cdot \mathcal{J}=\text{fibrations}$.
\end{itemize}}
\end{remark}

In theory of model categories, there exists an important notion of a $\mathcal{C}$-local object, which allows describing cases when the total left (resp. right) derived functor reverses arrows in the homotopy category (see Definition 3.1.4 \cite{Hr}). This definition goes back to the more general concept of localization and the related notion of a local object, which first appeared in \cite{GZ}. Let us recall it.

\begin{definition}\label{d6}
{\rm Let $\mathcal{C}$ be a category, and let $\mathcal{L}\subseteq\texttt{Mor}\ \mathcal{C}$ be a set of morphisms. An object $X \in \mathcal{C}$ is called a \textit{$\mathcal{L}$-local object} if for all $m \in \mathcal{L}$, the Hom functor from $m$ to $X$ induces a bijection
$\mathcal{C}(m, X)\colon\mathcal{C}(B, X) \xrightarrow{\cong}\mathcal{C}(A, X)$, that is, if every morphism $f\colon A \to X$ extends uniquely along $m$ to $B$:
$$
\xymatrix{
A \ar[r]^{f} \ar[d]_{m} & X \\
B \ar@{..>}[ur]_{\exists!} &
}
$$
We denote by \textit{$\mathcal{C}_{\mathcal{L}}$} the full \textit{subcategory of $\mathcal{L}$-local objects}}.
\end{definition}

\begin{remark}\label{r4}
{\rm If in Definition \ref{d6} the bijectivity condition is weakened to injectivity, then the object $X$ is called a \textit{$\mathcal{L}$-separating object}. The full \textit{subcategory of $\mathcal{L}$-separating objects} will be denoted by \textit{$\mathcal{C}^s_{\mathcal{L}}$}.}
\end{remark}

Let $(\mathcal{C}of,\mathcal{W},\mathcal{F}ib)$ be a qfs on category $\mathcal{C}$. Taking into account the adopted notation and by Definition \ref{d6}, we obtain $\mathcal{C}_f=\mathcal{C}_{\mathcal{C}of\cap\mathcal{W}}$.

A full subcategory $\mathcal{D}$ of category $\mathcal{C}$ is called \textit{reflective \emph{(}in $\mathcal{C}$\emph{)}} if the embedding functor
$\iota\colon\mathcal{D}\hookrightarrow\mathcal{C}$ has a left adjoint $L\colon\mathcal{C}\to\mathcal{D}$. In this case, the functor $L$ is called a \textit{reflector}, and the unit of adjunction $\eta\colon1_{\mathcal{C}}\Longrightarrow\iota L$ is called a \textit{reflection}. If the components of reflection $\eta_A$, $A\in\mathcal{C}$ belong to some class of morphisms $\mathcal{L}$ of category $\mathcal{C}$, closed under compositions with the class \texttt{Iso}, then the reflective subcategory $\mathcal{D}$ is called \textit{$\mathcal{L}$-reflective}. If $\mathcal{L}$ is the class $\texttt{Epi}$ (resp. $\texttt{Mono}$, $\texttt{Bim}$), then the reflective subcategory is called \textit{epireflective} (resp. \textit{monoreflective}, \textit{bireflective}).

\begin{example}\label{e1}
{\rm (i) The Stone-\v{C}ech compactification $\beta_X\colon X\to\iota\beta X$, $X\in\mathbf{Tych}$ is a dense embedding of $X$ into the compact space $\beta X$. This means that the category $\mathbf{CompHaus}$ is an $\mathcal{L}$-reflective subcategory in $\mathbf{Tych}$, where
      $\mathcal{L}$ is the class of dense embeddings in $\mathbf{Tych}$.

(ii) The category $\mathbf{Tych}$ is epireflective in $\mathbf{Top}$. The reflector $R\colon\mathbf{Top}\to\mathbf{Tych}$ is known as the Tychonoff functor.}
\end{example}

\begin{lemma}\label{l1}
Let $\mathcal{L}$ and $\mathcal{R}$ be arbitrary classes of morphisms in a category $\mathcal{C}$. Then the following conditions hold:
\begin{itemize}
  \item [(i)] a morphism with $\mathcal{L}$-local codomain is right orthogonal to the class $\mathcal{L}$ if and only if its domain is $\mathcal{L}$-local;
  \item [(ii)] if $\emph{\texttt{Iso}}\subseteq\mathcal{L}\cap\mathcal{R}$, then $\mathcal{C}_{\mathcal{R}\cdot\mathcal{L}}=\mathcal{C}_{\mathcal{L}}\cap\mathcal{C}_{\mathcal{R}}$.
\end{itemize}
If $(\mathcal{E},\mathcal{J},\mathcal{M})$ is a dfs in the category $\mathcal{C}$, then:
\begin{itemize}
  \item [(iii)] $(\mathcal{M}\cdot\mathcal{J})\cap\mathcal{C}_{\mathcal{J}}=\mathcal{M}\cap\mathcal{C}_{\mathcal{J}}$;
  \item [(iv)] the full subcategory of $\mathcal{J}\cdot\mathcal{E}$-local objects is reflective.
  \end{itemize}
In particular, if $(\mathcal{C}of, \mathcal{W},\mathcal{F}ib)$ is a qfs in $\mathcal{C}$, then
$\mathcal{F}ib\cap\mathcal{C}_{\mathcal{B}if}=\mathcal{F}ib\cap\mathcal{W}\cap\mathcal{C}_{\mathcal{B}if}$ and the full subcategory of $\mathcal{C}of$-local objects is reflective.
\end{lemma}

\begin{proof}
(i) Any diagram of the form $B \overset{\mathcal{L}\ni e}{\longleftarrow} A \overset{f}{\longrightarrow} C \overset{m\in\mathcal{L}^\perp}{\longrightarrow} D$ with $D$ being $\mathcal{L}$-local can be uniquely completed to a commutative diagram
$$
\xymatrix{
  A \ar[d]_{\mathcal{L}\ni e} \ar[r]^{f}        & C \ar[d]^{m\in\mathcal{L}^\perp}  \\
  B \ar@{.>}[ur]|-{\omega} \ar@{.>}[r]_{v} & D.             }
$$
This implies that $C$ is an $\mathcal{L}$-local object.

Conversely, if in a commutative square the object $C$ is $\mathcal{L}$-local, then there exists a unique arrow $\omega$ satisfying $\omega e=f$. The equality $v=m\omega$ follows from the $\mathcal{L}$-localness of $D$. Thus, $m\in\mathcal{L}^\perp$.

(ii) The relation $\mathcal{C}_{\mathcal{L}}\cap\mathcal{C}_{\mathcal{R}}\subseteq\mathcal{C}_{\mathcal{R}\cdot\mathcal{L}}$ follows from the commutative diagram
$$
\xymatrix@R=0.5cm{
A \ar[r]^(0.29999){f} \ar[d]_{\mathcal{L}\ni e}       &    D\in \mathcal{C}_{\mathcal{L}}\cap\mathcal{C}_{\mathcal{R}}         \\
B \ar@{.>}[ur]^(0.4){g} \ar[d]_{\mathcal{R}\ni m} &                                                                         \\
C \ar@{.>}[uur]_{h}      }
$$
From the condition $\texttt{Iso}\subseteq\mathcal{L}\cap\mathcal{R}$ we obtain $\mathcal{L}\subseteq\texttt{Iso}\cdot\mathcal{L}\subseteq\mathcal{R}\cdot\mathcal{L}$ and
$\mathcal{R}\subseteq\mathcal{R}\cdot\texttt{Iso}\subseteq\mathcal{R}\cdot\mathcal{L}$, hence $\mathcal{C}_{\mathcal{R}\cdot\mathcal{L}}\subseteq\mathcal{C}_{\mathcal{L}}\cap\mathcal{C}_{\mathcal{R}}$.

(iii) Let $m\colon A\to B$ be a fibration between $\mathcal{J}$-local objects. Consider the commutative rectangle
$$
\xymatrix{
  C \ar[d]_u \ar@<+0.5ex>[r]^{\mathcal{E}\ni e}         & E \ar@{.>}[dl]_k \ar@<+0.5ex>[r]^{j\in\mathcal{J}}  & D \ar@{.>}[dll]^h  \ar[d]^v \\
  A \ar[rr]_{m\in\mathcal{M}\cdot\mathcal{J}}                   && B.   }
$$
Since $\mathcal{M}\cdot\mathcal{J}=\mathcal{E}^{\bot}$, there exists a unique arrow $k$ satisfying $u=ke$ and $mk=vj$. By assumption $A\in\mathcal{C}_{\mathcal{J}}$, hence there exists a unique morphism $h$ such that $k=hj$. From the conditions $vj = mhj$, $j\in\mathcal{J}$ and $B\in\mathcal{C}_{\mathcal{J}}$, we obtain the equality $v=mh$. The uniqueness of $h$ follows from $m\in\mathcal{M}\cdot\mathcal{J}$ and $A\in\mathcal{C}_{\mathcal{J}}$.

(iv) Every object $X\in\mathcal{C}$ has an $(\mathcal{E},\mathcal{J},\mathcal{M})$-decomposition $$X \overset{e}{\to} A \overset{j}{\to} B \overset{!_B}{\to} 1.$$ Clearly, the terminal object 1 is $\mathcal{J}$-local. From part (i) and the relation $!_B\in\mathcal{M}=(\mathcal{J}\cdot\mathcal{E})^\perp$, it follows that $B$ is $\mathcal{J}\cdot\mathcal{E}$-local. The universality of the arrow $je\colon X\to B$ from $X$ to the embedding functor $\iota\colon\mathcal{C}_{\mathcal{J}\cdot\mathcal{E}}\hookrightarrow\mathcal{C}$ is obvious (see, for example, Proposition 5.5.5 part (i) \cite{B} or part 2.5 \cite{Bo}).
\end{proof}

\subsection{Double factorization systems and Quillen functors.}

In this section, we formulate the definition of Quillen functors, relying on similar definitions from model category theory. Detailed information about model categories and Quillen functors can be found in \cite{H}.

\begin{definition}\label{d7}
{\rm A\textit{model category $(\mathcal{C},\mathcal{C}of,\mathcal{W},\mathcal{F}ib)$} consists of a $\mathcal{C}$, together with three classes of morphisms $\mathcal{F}ib,\mathcal{C}of$ and $\mathcal{W}$, such that $\mathcal{W}$ satisfies the 2-out-of-3 property and $(\mathcal{C}of,\mathcal{F}ib\cap\mathcal{W})$ and $(\mathcal{C}of\cap\mathcal{W},\mathcal{F}ib)$ are weak factorization systems. If $(\mathcal{C}of,\mathcal{F}ib\cap\mathcal{W})$ and $(\mathcal{C}of\cap\mathcal{W},\mathcal{F}ib)$ are orthogonal factorization systems, then the model category $(\mathcal{C},\mathcal{C}of,\mathcal{W},\mathcal{F}ib)$ will be called \textit{orthogonal}}.
\end{definition}

We will use terminology for dfs analogous to that used in model category theory.

\begin{definition}\label{d8}
{\rm Let $(\mathcal{E}_{\mathcal{C}},\mathcal{J}_{\mathcal{C}},\mathcal{M}_{\mathcal{C}})$ and $(\mathcal{E}_{\mathcal{D}},\mathcal{J}_{\mathcal{D}},\mathcal{M}_{\mathcal{D}})$ be dfss on $\mathcal{C}$ and $\mathcal{D}$ respectively.
\begin{enumerate}
  \item [(i)] A functor $F\colon\mathcal{C}\to\mathcal{D}$ will be called a \textit{left Quillen functor} if $F$ is left adjoint and $F(\mathcal{J}_{\mathcal{C}}\cdot\mathcal{E}_{\mathcal{C}})\subseteq\mathcal{J}_{\mathcal{D}}\cdot\mathcal{E}_{\mathcal{D}}$ and $F(\mathcal{E}_{\mathcal{C}})\subseteq \mathcal{E}_{\mathcal{D}}$;
  \item [(ii)] A functor $G\colon\mathcal{D}\to\mathcal{C}$ will be called a \textit{right Quillen functor} if $G$ is right adjoint and
      $G(\mathcal{M}_{\mathcal{D}}\cdot\mathcal{J}_{\mathcal{D}})\subseteq\mathcal{M}_{\mathcal{C}}\cdot\mathcal{J}_{\mathcal{C}}$ and $G(\mathcal{M}_{\mathcal{D}})\subseteq\mathcal{M}_{\mathcal{C}}$;
  \item [(iii)] Let $F\colon\mathcal{C}\leftrightarrows\mathcal{D}\colon G$ be an adjunction from $\mathcal{C}$ to $\mathcal{D}$. We will call $F\colon\mathcal{C}\leftrightarrows\mathcal{D}\colon G$ a \textit{Quillen adjunction} if $F$ is a left Quillen functor.
\end{enumerate}}
\end{definition}

There is a statement analogous to Lemma 1.3.4. \cite{H}: $F\colon\mathcal{C}\leftrightarrows\mathcal{D}\colon G$ is a Quillen adjunction if and only if $G$ is a right Quillen functor.

\subsection{Right-induced double factorization systems on the category of algebras.}

We begin by recalling the notions of monads and $\mathbb{T}$-algebras (Eilenberg-Moore algebras). For more detailed information about monads and properties of the category of $\mathbb{T}$-algebras, the reader may found in any standard textbook on category theory (for example \cite{B}).

\begin{definition}\label{d9}
{\rm A \textit{monad} $\mathbb{T}$ in a category $\mathcal{C}$ is an ordered triple $\mathbb{T}=(T,\eta,\mu)$, consisting of a functor $T\colon\mathcal{C}\to\mathcal{C}$ and two natural transformations $\eta\colon 1_{\mathcal{C}}\Longrightarrow T$ and $\mu\colon T^2\Longrightarrow T$, called the \textit{unit} and \textit{multiplication} respectively, such that the following diagrams commute:}
$$
\xymatrix{T \ar@{=}[dr] \ar[r]^{T\eta} & T^2 \ar[d]^{\mu} & 1 \ar[l]_{\eta_T} \ar@{=}[dl]  \\
                & T,             }
\
\
\xymatrix{T^3 \ar[d]_{T\mu} \ar[r]^{\mu_T} & T^2 \ar[d]^{\mu}  \\
          T^2                    \ar[r]_{\mu} & T.             }
$$
\end{definition}

\begin{definition}\label{d10}
{\rm Let $\mathbb{T}$ be a monad on $\mathcal{C}$. An ordered pair $(A,h)$, consisting of an object $A$ and an arrow $h\colon TA\to A$, is called a \textit{$\mathbb{T}$-algebra} (\textit{Eilenberg-Moore algebra}) if the following diagram commutes:
$$
\xymatrix{A \ar@{=}[dr] \ar[r]^{\eta_A} & TA \ar[d]^{h} & T^2A \ar[l]_{Th} \ar[d]^{\mu_A}  \\
                                        & A             & TA   \ar[l]_{h}     .         }
$$
A \textit{morphism $f\colon(A,h)\to(B,v)$ of $\mathbb{T}$-algebras} from $(A,h)$ to $(B,v)$ is a morphism $f\colon A\to B$ in $\mathcal{C}$ such that the following diagram commutes:}
$$
\xymatrix{ TA \ar[d]_{h}  \ar[r]^{Tf}           & TB \ar[d]^{v}  \\
            A \ar[r]_{f}                        & B.             }
$$
\end{definition}
$\mathbb{T}$-algebras and morphisms of $\mathbb{T}$-algebras form a category $\mathcal{C}^\mathbb{T}$ (with composition and identity arrows defined in the obvious way).

We are interested in the following question. Under what conditions can a dfs on the category $\mathcal{C}$ be extended to the category of $\mathbb{T}$-algebras $\mathcal{C}^\mathbb{T}$ in such a way that $F^{\mathbb{T}}\colon\mathcal{C}\leftrightarrows\mathcal{C}^{\mathbb{T}}\colon V^{\mathbb{T}}$ forms a Quillen adjunction? It suffices to show that the forgetful functor $V^{\mathbb{T}}$ is a right Quillen functor. By analogy with the definition of induced model structure, we introduce the definition of an induced dfs.

\begin{definition}\label{d11}
{\rm Let $F\colon\mathcal{C}\rightleftarrows\mathcal{D}\colon G$ be an adjunction and $(\mathcal{E}_{\mathcal{C}},\mathcal{J}_{\mathcal{C}},\mathcal{M}_{\mathcal{C}})$ (resp. $(\mathcal{E}_{\mathcal{D}},\mathcal{J}_{\mathcal{D}},\mathcal{M}_{\mathcal{D}})$) be a dfs on $\mathcal{C}$ (resp. $\mathcal{D}$).
The system $({}^\perp(G^{-1}(\mathcal{M}\cdot\mathcal{J})),G^{-1}(\mathcal{J}),G^{-1}(\mathcal{M}))$ (resp.
$(F^{-1}(\mathcal{E}),F^{-1}(\mathcal{J}),(F^{-1}(\mathcal{J}\cdot\mathcal{E}))^\perp)$) is called the \textit{right-induced} (resp. \textit{left-induced}) \textit{dfs}, if it forms a dfs on $\mathcal{D}$ (resp. on $\mathcal{C}$).}
\end{definition}

Let $(\mathcal{E},\mathcal{J},\mathcal{M})$ be a dfs and $\mathbb{T}=(T,\eta,\mu)$ be a monad in $\mathcal{C}$. If $T$ preserves bifibrant morphisms and trivial cofibrations, then the dfs $(\mathcal{E},\mathcal{J},\mathcal{M})$ in $\mathcal{C}$ generates the right-induced dfs $(\mathcal{E}^{\mathbb{T}},\mathcal{J}^{\mathbb{T}},\mathcal{M}^{\mathbb{T}})$ in $\mathcal{C}^{\mathbb{T}}$. Let $\mathcal{E}^{\mathbb{T}}$, $\mathcal{J}^{\mathbb{T}}$ and $\mathcal{M}^{\mathbb{T}}$ be the preimages under the forgetful functor $V^{\mathbb{T}}$ of the classes $\mathcal{E}$, $\mathcal{J}$ and $\mathcal{M}$, respectively. We will show that the system defined in such a way forms a right-induced dfs in $\mathcal{C}^{\mathbb{T}}$.

\begin{lemma}\label{l2}
Let $\mathcal{L}$ and $\mathcal{R}$ be mutually orthogonal classes of morphisms and $\mathbb{T}=(T,\eta,\mu)$ be a monad in $\mathcal{C}$. Assume $\mathcal{L}^{\mathbb{T}}=(V^{\mathbb{T}})^{-1}(\mathcal{L})$, $\mathcal{R}^{\mathbb{T}}=(V^{\mathbb{T}})^{-1}(\mathcal{R})$ and $(\mathcal{R}\cdot\mathcal{L})^{\mathbb{T}}=(V^{\mathbb{T}})^{-1}(\mathcal{R}\cdot\mathcal{L})$. If $T$ satisfy conditions $T(\mathcal{L})\perp\mathcal{R}$ and $T^2(\mathcal{L})\perp\mathcal{R}$, in particular $T(\mathcal{L})\subseteq\mathcal{L}$, then $\mathcal{L}^{\mathbb{T}}\perp\mathcal{R}^{\mathbb{T}}$ and $(\mathcal{R}\cdot\mathcal{L})^{\mathbb{T}}=\mathcal{R}^{\mathbb{T}}\cdot\mathcal{L}^{\mathbb{T}}$. In addition, if $(\mathcal{L}, \mathcal{R})$ is an fs, then $(\mathcal{L}^{\mathbb{T}},\mathcal{R}^{\mathbb{T}})$ is also an fs.
\end{lemma}

\begin{proof}
The inclusion $\mathcal{R}^{\mathbb{T}}\cdot\mathcal{L}^{\mathbb{T}}\subseteq(\mathcal{R}\cdot\mathcal{L})^{\mathbb{T}}$ is obvious. Let us prove the reverse inclusion. Let $f\colon(A,h)\to(B,d)$ be a morphism of $\mathbb{T}$-algebras belonging to the class $(\mathcal{R}\cdot\mathcal{L})^{\mathbb{T}}$, and let $me$ be an $(\mathcal{L}, \mathcal{R})$-factorization of the morphism $f$. Since $Te\in T(\mathcal{L})\subseteq{}^\perp\mathcal{R}$ and $T^2e\in T^2(\mathcal{L})\subseteq{}^\perp\mathcal{R}$, there exists a unique arrow $k$ satisfying the conditions $mk=dTm$ and $kTe=eh$. From the equalities
\begin{align*}
   mk\mu_C     & = dTm\mu_C        &  k\mu_CT^2e & = kTe\mu_A       &    mk\eta_C   & = dTm\eta_C    &     k\eta_Ce    & = kTe\eta_A                 \\
               & = d\mu_BT^2m      &             & = eh\mu_A        &               & = d\eta_Bm     &                 & = eh\eta_A                    \\
               & = dTdT^2m         &             & = ehTh           &               & = m            &                 & = e                             \\
               & = dTmTk           &             & = kTeTh          &               &                &                 &                                      \\
               & = mkTk            &             & = kTkT^2e        &               &                &                 &                                         \\
\end{align*}
it follows that $(C,k)$ is a $\mathbb{T}$-algebra, and $me$ is an $(\mathcal{L}^{\mathbb{T}},\mathcal{R}^{\mathbb{T}})$-factorization of the morphism $f$ in $\mathcal{C}^{\mathbb{T}}$.
Therefore, $(\mathcal{R}\cdot\mathcal{L})^{\mathbb{T}}=\mathcal{R}^{\mathbb{T}}\cdot\mathcal{L}^{\mathbb{T}}$.

Let us prove the mutual orthogonality of the classes $\mathcal{L}^{\mathbb{T}}$ and $\mathcal{R}^{\mathbb{T}}$. From the commutative diagram
$$\xymatrix{
TA \ar[rr]^{Tu} \ar[dr]_{Te} \ar[dd]_s &                                    &          TC \ar[dr]^{Tm}  \ar'[d][dd]^(0.4)k      \\
                           &  TB  \ar[rr]^(0.34){Tv}  \ar[ur]^{T\omega} \ar[dd]_(0.29)t &       &  TD  \ar[dd]^d                   \\
A \ar'[r][rr]^(0.4)u \ar[dr]_{\mathcal{L}^{\mathbb{T}}\ni e}      &                                    &          C    \ar[dr]^(0.4){m\in\mathcal{R}^{\mathbb{T}}}                   \\
                           &  B  \ar@{.>}[ur]_{\omega}  \ar[rr]_v         &       &  D                               \\
}
$$
it follows that there exists a unique morphism $\omega\colon(B,t)\to(C,k)$ satisfying the equalities $\omega e=u$ and $m\omega=v$ in $\mathcal{C}^{\mathbb{T}}$. Thus, $\mathcal{L}^{\mathbb{T}}\perp\mathcal{R}^{\mathbb{T}}$.

In addition, if $(\mathcal{L}, \mathcal{R})$ is an fs, then $(\mathcal{L}^{\mathbb{T}},\mathcal{R}^{\mathbb{T}})$ satisfy conditions (i)-(iv) of Definition 5.5.1 Vol.1 \cite{B}. By Proposition 5.5.3 Vol.1 \cite{B} $(\mathcal{L}^{\mathbb{T}},\mathcal{R}^{\mathbb{T}})$ is an fs.
\end{proof}

\begin{proposition}\label{p1}
Let $(\mathcal{E}, \mathcal{J}, \mathcal{M})$ be a dfs in a category $\mathcal{C}$, $\mathbb{T}=(T,\eta,\mu)$ be a monad in $\mathcal{C}$, and the functor $T$ preserve all small \emph{(}finite\emph{)} direct limits, trivial cofibrations, and bifibrant morphisms. Let $\mathcal{S}$ denote one of the classes $\mathcal{J}$ or $\mathcal{J}\cdot\mathcal{E}$. Then the following conditions hold:
\begin{itemize}
  \item [(i)] $V^{\mathbb{T}}$ generates the right-induced dfs $(\mathcal{E}^{\mathbb{T}},\mathcal{J}^{\mathbb{T}},\mathcal{M}^{\mathbb{T}})$ in $\mathcal{C}^{\mathbb{T}}$;
  \item [(ii)] $V^{\mathbb{T}}$ preserves and reflects the corresponding local objects, i.e., for any object $(A,h)\in C^{\mathbb{T}}$, the following equivalence holds\emph{:} $$(A,h)\in\mathcal{C}^{\mathbb{T}}_{\mathcal{S}^{\mathbb{T}}}\Longleftrightarrow A\in\mathcal{C}_{\mathcal{S}}.$$
\end{itemize}
In particular, $F^{\mathbb{T}}\colon\mathcal{C}\leftrightarrows\mathcal{C}^{\mathbb{T}}\colon V^{\mathbb{T}}$ is a Quillen adjunction.
\end{proposition}

\begin{proof}
(i) Note that the (finite) bicompleteness of the category of $\mathbb{T}$-algebras is guaranteed by the preservation of all small (finite) direct limits by $T$ (see Proposition 4.3.1 and Proposition 4.3.2 \cite{B}). By assumption, $T$ preserves trivial cofibrations and bifibrant morphisms. Consequently, $T$ preserves cofibrations. Define $(\mathcal{E}^{\mathbb{T}},\mathcal{J}^{\mathbb{T}},\mathcal{M}^{\mathbb{T}})$ as follows: $$\mathcal{E}^{\mathbb{T}}=(V^{\mathbb{T}})^{-1}(\mathcal{E}),\ \mathcal{J}^{\mathbb{T}}=(V^{\mathbb{T}})^{-1}(\mathcal{J})\ \text{and}\ \mathcal{M}^{\mathbb{T}}=(V^{\mathbb{T}})^{-1}(\mathcal{M}).$$ By Proposition 2.3 \cite{PT}, $(\mathcal{E},\mathcal{M}\cdot\mathcal{J})$ and $(\mathcal{J}\cdot\mathcal{E},\mathcal{M})$ are fs's and $\mathcal{E}\perp\mathcal{M}$. Since $\mathcal{J}={}^\perp\mathcal{M}\cap\mathcal{E}^\perp$, we have $\mathcal{E}\perp\mathcal{J}$ and $\mathcal{J}\perp\mathcal{M}$. From Lemma \ref{l2}, it follows that $(\mathcal{J}\cdot\mathcal{E})^{\mathbb{T}}=\mathcal{J}^{\mathbb{T}}\cdot\mathcal{E}^{\mathbb{T}}$,
$(\mathcal{M}\cdot\mathcal{J})^{\mathbb{T}}=\mathcal{M}^{\mathbb{T}}\cdot\mathcal{J}^{\mathbb{T}}$ and $\mathcal{E}^{\mathbb{T}}\perp\mathcal{M}^{\mathbb{T}}$, thus $(\mathcal{E}^{\mathbb{T}},\mathcal{M}^{\mathbb{T}}\cdot\mathcal{J}^{\mathbb{T}})$ and $(\mathcal{J}^{\mathbb{T}}\cdot\mathcal{E}^{\mathbb{T}},\mathcal{M}^{\mathbb{T}})$ are fs's. Since $(\mathcal{M}^{\mathbb{T}}\cdot\mathcal{J}^{\mathbb{T}})\cap (\mathcal{J}^{\mathbb{T}}\cdot\mathcal{E}^{\mathbb{T}})=\mathcal{J}^{\mathbb{T}}$, by Theorem 2.7 \cite{PT}, $(\mathcal{E}^{\mathbb{T}},\mathcal{J}^{\mathbb{T}},\mathcal{M}^{\mathbb{T}})$ is a right-induced dfs.

(ii) Let $V^{\mathbb{T}}(A,h)=A\in\mathcal{C}_{\mathcal{S}}$. Then for any arrow $f\colon(B,k)\to(A,h)$ and any cofibration or bifibrant morphism $e\colon(B,k)\to(C,u)$ (since $V^{\mathbb{T}}$ generates the dfs $(\mathcal{E}^{\mathbb{T}},\mathcal{J}^{\mathbb{T}},\mathcal{M}^{\mathbb{T}})$, we have $e\colon(B,k)\to(C,u)\in\mathcal{S}^{\mathbb{T}}\Longrightarrow e\in\mathcal{S}\Longrightarrow Te\in\mathcal{S}\Longrightarrow Te\colon(TB,\mu_B)\to(TC,\mu_C)\in\mathcal{S}^{\mathbb{T}}$), there exists a unique arrow \(g\colon C\to A\) such that $f=ge$ in $\mathcal{C}$. From the equality $guTe=gek=fk=hTf=hTgTe$, it follows that $g\colon(C,u)\to(A,h)$ is a morphism of $\mathbb{T}$-algebras. Thus, $V^{\mathbb{T}}$ reflects $\mathcal{S}$-local objects.

Now suppose $(A,h)$ is an $\mathcal{S}^{\mathbb{T}}$-local object in $\mathcal{C}^{\mathbb{T}}$. For any $f\colon B\to A$ and any cofibration or bifibrant morphism $e\colon B\to C$, there exists a unique arrow $g\colon(TC,\mu_C)\to(A,h)$ making the diagrams
$$
\xymatrix{
  (TB,\mu_B) \ar[dd]_{Tf} \ar[rr]^{Te}      & & (TC,\mu_C)  \ar[ddll]<+0.5ex>^{T\omega_1}\ar[ddll]<-0.5ex>_{T\omega_2}\ar@{.>}[dd]^{g} \\
                                          &&                                                                                         \\
  (TA,\mu_A) \ar[rr]_h & &(A,h)   }
\
\
\
\
\xymatrix{
B \ar[rr]^e \ar[dr]_f \ar[dd]_{\eta_B} &                                    &          C \ar[dl]<+0.5ex>^{\omega_1}\ar[dl]<-0.5ex>_{\omega_2} \ar[dr]^{g\eta_C}  \ar'[d][dd]^(0.4){\eta_C}      \\
                                          &  A  \ar@{=}[rr]   \ar[dd]_(0.29){\eta_A} &       &  A  \ar@{=}[dd]                   \\
TB\ar'[r][rr]^(0.4){Te} \ar[dr]_{Tf}      &                                    &          TC  \ar[dl]<+0.5ex>^{T\omega_1}\ar[dl]<-0.5ex>_{T\omega_2}  \ar[dr]^(0.4){g}                  \\
                           &  TA    \ar[rr]_h         &       &  A                               \\
}
$$
commutative. From these diagrams, it follows that $g\eta_C$ is the unique arrow in $\mathcal{C}$ satisfying $g\eta_Ce=f$. Therefore, $V^{\mathbb{T}}$ maps $\mathcal{S}^{\mathbb{T}}$-local objects to $\mathcal{S}$-local objects.
\end{proof}

\begin{remark}\label{r5}
{\rm In the proof of Proposition \ref{p1}, the preservation of small (finite) direct limits by $T$ was used to justify the (finite) cocompleteness of the category $\mathcal{C}^{\mathbb{T}}$. Proposition \ref{p1} remains valid if the preservation of small (finite) direct limits by $T$ is replaced by (finite) cocompletness of the category $\mathcal{C}^{\mathbb{T}}$.}
\end{remark}

\begin{corollary}\label{c1}
Let $(\mathcal{C},\mathcal{C}of,\mathcal{W},\mathcal{F}ib)$ be an orthogonal model category, \(\mathbb{T}=(T,\eta,\mu)\) be a monad in \(\mathcal{C}\), and the functor \(T\) preserve all small \emph{(}finite\emph{)} direct limits, trivial cofibrations and bifibrant morphisms. Then the following conditions hold:
\begin{itemize}
  \item [(i)] $V^{\mathbb{T}}$ generates a right-induced qfs $(\mathcal{C}of^{\mathbb{T}},\mathcal{W}^{\mathbb{T}},\mathcal{F}ib^{\mathbb{T}})$;
  \item [(ii)] $V^{\mathbb{T}}$ preserves and reflects the corresponding $\mathcal{C}of$-local objects, i.e., for any object $(A,h)\in V^{\mathbb{T}}$, the following equivalence holds\emph{:} $$(A,h)\in\mathcal{C}^{\mathbb{T}}_{\mathcal{C}of^{\mathbb{T}}}\Longleftrightarrow A\in\mathcal{C}_{\mathcal{C}of}.$$
\end{itemize}
In particular, $(\mathcal{C}^{\mathbb{T}},\mathcal{C}of^{\mathbb{T}},\mathcal{W}^{\mathbb{T}},\mathcal{F}ib^{\mathbb{T}})$ is a orthogonal model category, and $F^{\mathbb{T}}\colon\mathcal{C}\leftrightarrows\mathcal{C}^{\mathbb{T}}\colon V^{\mathbb{T}}$ is a Quillen adjunction.
\end{corollary}

\begin{proof}
Define $(\mathcal{C}of^{\mathbb{T}},\mathcal{W}^{\mathbb{T}},\mathcal{F}ib^{\mathbb{T}})$ as follows: $$\mathcal{W}^{\mathbb{T}}=(V^{\mathbb{T}})^{-1}(\mathcal{W}),\ \mathcal{F}ib^{\mathbb{T}}=(V^{\mathbb{T}})^{-1}(\mathcal{F}ib)\ \text{and}\ \mathcal{C}of^{\mathbb{T}}=(V^{\mathbb{T}})^{-1}(\mathcal{C}of).$$
The 2-out-of-3 property follows directly from the definition of the class \(\mathcal{W}^{\mathbb{T}}\). By Theorem \ref{th1}, the qfs $(\mathcal{C}of,\mathcal{W},\mathcal{F}ib)$ corresponds to the dfs $(\mathcal{C}of\cap\mathcal{W},\mathcal{F}ib\cap\mathcal{C}of,\mathcal{F}ib\cap\mathcal{W})$.
Since $\mathcal{C}of=(\mathcal{F}ib\cap\mathcal{C}of)\cdot(\mathcal{C}of\cap\mathcal{W})$, the functor $T$ preserves the class $\mathcal{C}of$. By Proposition \ref{p1}, $((\mathcal{C}of\cap\mathcal{W})^{\mathbb{T}},(\mathcal{F}ib\cap\mathcal{C}of)^{\mathbb{T}},(\mathcal{F}ib\cap\mathcal{W})^{\mathbb{T}})$ is a right-induced dfs, where $(\mathcal{C}of\cap\mathcal{W})^{\mathbb{T}}=\mathcal{C}of^{\mathbb{T}}\cap\mathcal{W}^{\mathbb{T}}$,
$(\mathcal{F}ib\cap\mathcal{C}of)^{\mathbb{T}}=(\mathcal{F}ib)^{\mathbb{T}}\cap(\mathcal{C}of)^{\mathbb{T}}$,
$(\mathcal{F}ib\cap\mathcal{W})^{\mathbb{T}}=\mathcal{F}ib^{\mathbb{T}}\cap\mathcal{W}^{\mathbb{T}}$.
Thus, $(\mathcal{C}of^{\mathbb{T}},\mathcal{W}^{\mathbb{T}},\mathcal{F}ib^{\mathbb{T}})$ is a right-induced qfs (by Theorem 2.7 \cite{PT}, $(\mathcal{C}of^{\mathbb{T}}\cap\mathcal{W}^{\mathbb{T}},\mathcal{F}ib^{\mathbb{T}})$ and $(\mathcal{C}of^{\mathbb{T}},\mathcal{F}ib^{\mathbb{T}}\cap\mathcal{W}^{\mathbb{T}})$ are fs's), corresponding to the dfs
$(\mathcal{C}of^{\mathbb{T}}\cap\mathcal{W}^{\mathbb{T}},\mathcal{F}ib^{\mathbb{T}}\cap\mathcal{C}of^{\mathbb{T}},
\mathcal{F}ib^{\mathbb{T}}\cap\mathcal{W}^{\mathbb{T}})$. In particular, $(\mathcal{C}^{\mathbb{T}},\mathcal{C}of^{\mathbb{T}},\mathcal{W}^{\mathbb{T}},\mathcal{F}ib^{\mathbb{T}})$ is a Quillen model category, and $F^{\mathbb{T}}\colon\mathcal{C}\leftrightarrows\mathcal{C}^{\mathbb{T}}\colon V^{\mathbb{T}}$ is a Quillen adjunction. To complete the proof, it remains to refer to Proposition \ref{p1}.
\end{proof}

\begin{corollary}\label{c2}
Let $(\mathcal{E}_{\mathcal{C}},\mathcal{J}_{\mathcal{C}},\mathcal{M}_{\mathcal{C}})$ and $(\mathcal{E}_{\mathcal{D}},\mathcal{J}_{\mathcal{D}},\mathcal{M}_{\mathcal{D}})$ be dfs's, $\mathbb{T}=(T,\eta,\mu)$ and $\mathbb{H}=(H,\zeta,\nu)$ be monads in $\mathcal{C}$ and $\mathcal{D}$ respectively. Let the functors $T$ and $H$ preserve all small \emph{(}finite\emph{)} direct limits, trivial cofibrations and bifibrant morphisms. If in the commutative diagram
$$
\xymatrix{\mathcal C^{\mathbb T}\ar@<-0.5ex>[d]_{V^{\mathbb{T}}}\ar@<-0.5ex>[r]_P &  \mathcal D^{\mathbb H}\ar@<-0.5ex>[d]_{V^{\mathbb{H}}}\ar@<-0.5ex>[l]_Q    \\
              \mathcal C\ar@<-0.5ex>[r]_L\ar@<-0.5ex>[u]_{F^{\mathbb{T}}}           &   \mathcal D\ar@<-0.5ex>[u]_{F^{\mathbb{H}}}\ar@<-0.5ex>[l]_R          }
$$
$L\colon\mathcal{C}\leftrightarrows\mathcal{D}\colon R$ is a Quillen adjunction, then $P\colon\mathcal{C}^{\mathbb{T}}\leftrightarrows\mathcal{D}^{\mathbb{H}}\colon Q$ is also a Quillen adjunction.
\end{corollary}

\begin{proof} We notice that the equality $RV^{\mathbb{H}}=V^{\mathbb{T}}Q$ guarantees the existence of a left adjoint $P$ to the functor $Q$ and the commutativity (up to natural isomorphism) of the left adjoints. Indeed, since coequalizers exist in the category $\mathcal{D}^{\mathbb{H}}$, by Theorem 4.5.6. Vol.2 \cite{B}, the existence of a left adjoint $P$ to the functor $Q$ follows from the existence of a left adjoint $L$ to the functor $R$. It remains to refer to Proposition \ref{p1} item (ii).
\end{proof}

\begin{remark}\label{r6}
{\rm The diagram in Corollary \ref{c2} commutes with respect to right adjoint functors. This implies commutativity (up to natural isomorphism) with respect to left adjoint functors. The adjoint $P\colon\mathcal{C}^{\mathbb{T}}\leftrightarrows\mathcal{D}^{\mathbb{H}}\colon Q$ will be called the \textit{extension of the adjoint} $L\colon\mathcal{C}\leftrightarrows\mathcal{D}\colon R$.}
\end{remark}

\subsection{Left-induced double factorization systems on the category of coalgebras.}

The dual notion to an algebra is that of a coalgebra. In other words, if we reverse all arrows in the commutative diagrams in Definition \ref{d10} and Definition \ref{d11}, we obtain the definitions of comonad and coalgebra respectively. For more detailed information, the reader may refer to standard textbooks on category theory, such as \cite{MM} or \cite{M}.

By the principle of duality, the results obtained above automatically carry over to the case of the category of coalgebras. In this case, $\mathcal{J}\cdot\mathcal{E}$-local objects should be replaced with $\mathcal{M}\cdot\mathcal{J}$-colocal ones. However, of greater interest is studying the relationship between the category of $\mathcal{J}\cdot\mathcal{E}$-local objects and the category of coalgebras. We begin with the following auxiliary lemma.

\begin{lemma}\label{l3}
Let $(\mathcal{E},\mathcal{J},\mathcal{M})$ be a dfs in the category $\mathcal{C}$. If $A$ is a $\mathcal{J}\cdot\mathcal{E}$-local object, then the diagonal monomorphism $\Delta_A=(1_A,1_A)\colon A\rightarrowtail A\times A$ is a trivial fibration. Conversely, if $\Delta_A$ is a trivial fibration, then $A$ is a $\mathcal{J}\cdot\mathcal{E}$-separating object.
\end{lemma}

\begin{proof}
Let $A$ be a $\mathcal{J}\cdot\mathcal{E}$-local object. Then $A\times A$ is also a $\mathcal{J}\cdot\mathcal{E}$-local object (direct verification). Consider the $(\mathcal{E},\mathcal{J},\mathcal{M})$-factorization $A \overset{e}{\to} B \overset{j}{\to} C \overset{m}{\to} A\times A$ of the arrow $\Delta_A$. By Lemma \ref{l1} item (i), we obtain that $C$ is a $\mathcal{J}\cdot\mathcal{E}$-local object. From Lemma \ref{l1} item (iv) it follows that $je\colon A\to C$ is the universal arrow from $A$ to the embedding functor $\iota\colon\mathcal{C}_{\mathcal{J}\cdot\mathcal{E}}\hookrightarrow\mathcal{C}$. Since $\mathcal{C}_{\mathcal{J}\cdot\mathcal{E}}$ is a full reflective subcategory of $\mathcal{C}$, the unit of the adjunction on objects from $\mathcal{C}_{\mathcal{J}\cdot\mathcal{E}}$ is an isomorphism, hence $je\colon A\simeq C$ is an isomorphism in the category $\mathcal{C}$. Thus, $\Delta_A$ is a trivial fibration.

Conversely, let $\Delta_A$ be a trivial fibration. Consider arbitrary morphisms $f,g\colon B\to A$ and a cofibration $je\colon C\to B$ such that $fje=gje$. The required equality $f=g$ follows from the commutative diagram
$$
\xymatrix{
  C \ar[d]_{\mathcal{J}\cdot\mathcal{E}\ni je} \ar[r]^{fje=gje} &  A \ar[d]^{\Delta_A\in\mathcal{M}={}^\perp(\mathcal{J}\cdot\mathcal{E})}  \\
  B \ar@{.>}[ur]|-{\omega} \ar[r]_{(f,g)}      &  A\times A      .       }
$$
\end{proof}

The following theorem establishes the relationship between $\mathcal{J}\cdot\mathcal{E}$-local objects and the category of $\mathbb{G}$-coalgebras $\mathcal{C}_{\mathbb{G}}$.

\begin{proposition}\label{p2}
Let $(\mathcal{E},\mathcal{J},\mathcal{M})$ be a dfs in the category $\mathcal{C}$, $\mathbb{G}=(G,\epsilon,\delta)$ be a comonad in $\mathcal{C}$, and the functor $G$ preserve all small \emph{(}finite\emph{)} limits, trivial fibrations and bifibrant morphisms. Then the following conditions hold\emph{:}
\begin{itemize}
  \item [(i)] $G_{\mathbb{G}}$ generates a left-induced dfs $(\mathcal{E}_{\mathbb{G}},\mathcal{J}_{\mathbb{G}},\mathcal{M}_{\mathbb{G}})$\emph{;}
  \item [(ii)] $G_{\mathbb{G}}$ reflects the corresponding $\mathcal{J}\cdot\mathcal{E}$-local objects, i.e. $G_{\mathbb{G}}^{-1}(\mathcal{C}_{\mathcal{J}\cdot\mathcal{E}})\subseteq(\mathcal{C}_{\mathbb{G}})_{\mathcal{J}_{\mathbb{G}}\cdot\mathcal{E}_{\mathbb{G}}}$\emph{;}
  \item [(iii)] $G_{\mathbb{G}}$ and $R_{\mathbb{G}}$ preserve the corresponding local objects if and only if $G$ preserves $\mathcal{J}\cdot\mathcal{E}$-local objects, i.e. the following equivalence holds\emph{:} $$G_{\mathbb{G}}((\mathcal{C}_{\mathbb{G}})_{\mathcal{J}_{\mathbb{G}}\cdot\mathcal{E}_{\mathbb{G}}})\subseteq\mathcal{C}_{\mathcal{J}\cdot\mathcal{E}}\ \&\
      R_{\mathbb{G}}(\mathcal{C}_{\mathcal{J}\cdot\mathcal{E}})\subseteq(\mathcal{C}_{\mathbb{G}})_{\mathcal{J}_{\mathbb{G}}\cdot\mathcal{E}_{\mathbb{G}}}\Longleftrightarrow G(\mathcal{C}_{\mathcal{J}\cdot\mathcal{E}})\subseteq\mathcal{C}_{\mathcal{J}\cdot\mathcal{E}}.$$
\end{itemize}
In particular, $G_{\mathbb{G}}\colon\mathcal{C}_{\mathbb{G}}\leftrightarrows\mathcal{C}\colon R_{\mathbb{G}}$ is a Quillen adjunction.
\end{proposition}

\begin{proof}
(i) By assumption, the functor $G$ preserves small (finite) limits, trivial fibrations and bifibrant morphisms, hence $G$ preserves fibrations. Define $(\mathcal{E}_{\mathbb{G}},\mathcal{J}_{\mathbb{G}},\mathcal{M}_{\mathbb{G}})$ as follows: $$\mathcal{E}_{\mathbb{G}}=(G_{\mathbb{G}})^{-1}(\mathcal{E}),\ \mathcal{J}_{\mathbb{G}}=(G_{\mathbb{G}})^{-1}(\mathcal{J})\ \text{and}\ \mathcal{M}_{\mathbb{G}}=(G_{\mathbb{G}})^{-1}(\mathcal{M}).$$
By Proposition 2.3 \cite{PT}, $(\mathcal{E},\mathcal{M}\cdot\mathcal{J})$ and $(\mathcal{J}\cdot\mathcal{E},\mathcal{M})$ are fs's and $\mathcal{E}\perp\mathcal{M}$. Since $\mathcal{J}={}^\perp\mathcal{M}\cap\mathcal{E}^\perp$, we have $\mathcal{E}\perp\mathcal{J}$ and $\mathcal{J}\perp\mathcal{M}$. The dual statement of Lemma \ref{l2} implies the equalities $(\mathcal{J}\cdot\mathcal{E})_{\mathbb{G}}=\mathcal{J}_{\mathbb{G}}\cdot\mathcal{E}_{\mathbb{G}}$,
$(\mathcal{M}\cdot\mathcal{J})_{\mathbb{G}}=\mathcal{M}_{\mathbb{G}}\cdot\mathcal{J}_{\mathbb{G}}$ and $\mathcal{E}_{\mathbb{G}}\perp\mathcal{M}_{\mathbb{G}}$, and hence $(\mathcal{E}_{\mathbb{G}},\mathcal{M}_{\mathbb{G}}\cdot\mathcal{J}_{\mathbb{G}})$ and $(\mathcal{J}_{\mathbb{G}}\cdot\mathcal{E}_{\mathbb{G}},\mathcal{M}_{\mathbb{G}})$ are fs's. Since $(\mathcal{M}_{\mathbb{G}}\cdot\mathcal{J}_{\mathbb{G}})\cap (\mathcal{J}_{\mathbb{G}}\cdot\mathcal{E}_{\mathbb{G}})=\mathcal{J}_{\mathbb{G}}$, then by Theorem 2.7 \cite{PT}, $(\mathcal{E}_{\mathbb{G}},\mathcal{J}_{\mathbb{G}},\mathcal{M}_{\mathbb{G}})$ is a left-induced dfs.

(ii) Let $G_{\mathbb{G}}(A,s)=A\in\mathcal{C}_{\mathcal{J}\cdot\mathcal{E}}$. Then for an arbitrary arrow $f\colon(B,t)\to(A,s)$ and any cofibration $v\colon(B,t)\to(C,u)$, there exists a unique arrow $g\colon C\to A$ such that $f=ge$ in $\mathcal{C}$. Since $A$ is a $\mathcal{J}\cdot\mathcal{E}$-local object, by Lemma \ref{l1} item (i), it follows that $!_A\colon A\to 1$ is a trivial fibration. By assumption, $G$ preserves terminal objects and trivial fibrations. Therefore, $G(!_A)=!_{GA}$ is a trivial fibration, hence by Lemma \ref{l1} item (i), $GA$ is a $\mathcal{J}\cdot\mathcal{E}$-local object. From the equality $G(g)ue=G(g)G(e)t=G(f)t=sf=sge$ we obtain that $g\colon(C,u)\to(A,s)$ is a morphism of $\mathbb{G}$-coalgebras. Thus, $G_{\mathbb{G}}$ reflects the corresponding $\mathcal{J}\cdot\mathcal{E}$-local objects.

(iii) Suppose the condition $G(\mathcal{C}_{\mathcal{J}\cdot\mathcal{E}})\subseteq\mathcal{C}_{\mathcal{J}\cdot\mathcal{E}}$ holds. Consider an arbitrary $\mathcal{J}_{\mathbb{G}}\cdot\mathcal{E}_{\mathbb{G}}$-local object $(A,s)$ in $\mathcal{C}_{\mathbb{G}}$. Let $\eta\colon1_{\mathcal{C}}\Longrightarrow\iota L$ be the unit of the adjunction $L\colon\mathcal{C}\leftrightarrows\mathcal{C}_{\mathcal{J}\cdot\mathcal{E}}\colon\iota$, and $\eta_{\mathbb{G}}\colon1_{\mathcal{C}_{\mathbb{G}}}\Longrightarrow\iota_{\mathbb{G}}L_{\mathbb{G}}$ be the unit of the adjunction $L_{\mathbb{G}}\colon\mathcal{C}_{\mathbb{G}}\leftrightarrows(\mathcal{C}_{\mathbb{G}})_{\mathcal{J}_{\mathbb{G}}\cdot\mathcal{E}_{\mathbb{G}}}\colon\iota_{\mathbb{G}}$. Since $\mathcal{C}_{\mathcal{J}\cdot\mathcal{E}}\ni\iota LA$, then $G\iota LA\in\mathcal{C}_{\mathcal{J}\cdot\mathcal{E}}$. Consequently, there exists a unique arrow $t$ such that $G(\eta_A)s=t\eta_A$. It is easy to verify that $(\iota LA,t)$ is a $\mathbb{G}$-coalgebra. For example, the equality $\delta_{\iota LA}t=G(t)t$ follows from the diagram
$$
\xymatrix{
  A \ar[dd]_{s} \ar[rd]^s \ar[rr]^{\eta_A} && \iota LA \ar[dd]^(0.666)t \ar[rd]^t  &            \\
  & GA \ar[dd]_(0.299999){\delta_A} \ar[rr]^(0.28){G\eta_A} & & G\iota LA\ar[dd]^{\delta_{\iota LA}}                  \\
  GA   \ar[rd]_{Gs} \ar[rr]^(0.666){G\eta_A}         && G\iota LA \ar[rd]^{Gt}  &                                   \\
  & G^2A  \ar[rr]_{G^2\eta_A}         && G^2\iota LA              }
$$
and the fact that $\eta$ is the unit of the adjunction $L\colon\mathcal{C}\leftrightarrows\mathcal{C}_{\mathcal{J}\cdot\mathcal{E}}\colon\iota$
and $G^2\iota LA\in\mathcal{C}_{\mathcal{J}\cdot\mathcal{E}}$. By item (ii), we obtain $(\iota LA,t)\in(\mathcal{C}_{\mathbb{G}})_{\mathcal{J}_{\mathbb{G}}\cdot\mathcal{E}_{\mathbb{G}}}$. It follows that $\eta_A\colon(A,s)\to(\iota LA,t)$ is a morphism in the category of $\mathbb{G}$-coalgebras.

Let $mje$ be the $(\mathcal{E}_{\mathbb{G}},\mathcal{J}_{\mathbb{G}},\mathcal{M}_{\mathbb{G}})$-factorization of the morphism $\eta_A$ in $\mathcal{C}_{\mathbb{G}}$. Since the functor $G_{\mathbb{G}}$ maps $(\mathcal{E}_{\mathbb{G}},\mathcal{J}_{\mathbb{G}},\mathcal{M}_{\mathbb{G}})$-factorizations to $(\mathcal{E},\mathcal{J},\mathcal{M})$-factorizations (follows from item (i)), then $mje$ is an $(\mathcal{E},\mathcal{J},\mathcal{M})$-factorization of the morphism $\eta_A$ in $\mathcal{C}$. From Lemma \ref{l1} item (iv) it follows that the unit of the adjunction $\eta_A$ has an $(\mathcal{E},\mathcal{J})$-factorization, hence $m$ is an isomorphism in $\mathcal{C}$. Thus, the morphism $\eta_A$ has an $(\mathcal{E}_{\mathbb{G}},\mathcal{J}_{\mathbb{G}})$-factorization in $\mathcal{C}_{\mathbb{G}}$. By Lemma \ref{l1} item (iv), we obtain that $\eta_A$ is the universal arrow from $(A,s)$ to the embedding functor $\iota_{\mathbb{G}}$. The latter means that there exists an isomorphism $\xi$ such that $\eta_A=\xi\eta_{\mathbb{G}}{}_{(A,s)}$. Since the category $(\mathcal{C}_{\mathbb{G}})_{\mathcal{J}_{\mathbb{G}}\cdot\mathcal{E}_{\mathbb{G}}}$ is a full reflective subcategory of $\mathcal{C}_{\mathbb{G}}$ and $(A,s)\in(\mathcal{C}_{\mathbb{G}})_{\mathcal{J}_{\mathbb{G}}\cdot\mathcal{E}_{\mathbb{G}}}$, the component $\eta_{\mathbb{G}}{}_{(A,s)}\colon(A,s)\simeq\iota_{\mathbb{G}}L_{\mathbb{G}}(A,s)$ of the unit of the adjunction $L_{\mathbb{G}}\colon\mathcal{C}_{\mathbb{G}}\leftrightarrows(\mathcal{C}_f\cap\mathcal{C}_{\mathcal{B}if})_{\mathbb{G}}\colon\iota_{\mathbb{G}}$ is an isomorphism. It follows that $\eta_A\colon(A,s)\simeq(\iota LA,t)$ is an isomorphism in $\mathcal{C}_{\mathbb{G}}$, and hence in $\mathcal{C}$. Finally, we obtain that $G_{\mathbb{G}}$ maps $\mathcal{J}_{\mathbb{G}}\cdot\mathcal{E}_{\mathbb{G}}$-local objects to $\mathcal{J}\cdot\mathcal{E}$-local ones. Since $G=G_{\mathbb{G}}R_{\mathbb{G}}$ and $G_{\mathbb{G}}$ reflects the corresponding $\mathcal{J}\cdot\mathcal{E}$-local objects, then $R_{\mathbb{G}}(\mathcal{C}_{\mathcal{J}\cdot\mathcal{E}})\subseteq(\mathcal{C}_{\mathbb{G}})_{\mathcal{J}_{\mathbb{G}}\cdot\mathcal{E}_{\mathbb{G}}}$.
\end{proof}

\section{Applications in topos theory}

\subsection{Basic concepts and definitions}

Here we recall the basic concepts and definitions from topos theory. For more detailed information about toposes, the reader may refer to \cite{MM}, \cite{J}.

\begin{definition}[Definition IV.1.2]\label{d12}
{\rm A \textit{topos} $\mathcal{C}$ is a category satisfying the following axioms:
\begin{itemize}
  \item [(i)] $\mathcal{C}$ has pullbacks;
  \item [(ii)] $\mathcal{C}$ has a terminal object 1;
  \item [(iii)] $\mathcal{C}$ has a \textit{subobject classifier}, i.e., there is a monomorphism $\text{true}\colon 1\rightarrowtail\Omega$ satisfying the following universal property: for any monomorphism \(m\colon A\rightarrowtail B\) there exists a unique arrow $\text{char}\ m$ making the diagram
      $$\xymatrix{ A \ar@{>->}[d]_{m} \ar[r]^{!_A}   & 1 \ar@{>->}[d]^{\text{true}}  \\
                   B  \ar@{.>}[r]_{\text{char}\ m}  & \Omega             }
   $$
      a pullback;
  \item [(iv)] $\mathcal{C}$ has exponential objects, i.e., for any object $A$ there exists a morphism $\epsilon_A\colon A\times PA\to\Omega$ satisfying the following universal property: for each morphism $f\colon A\times B\to\Omega$ there exists a unique arrow $g$ making the diagram
      $$\xymatrix{
   A\times PA  \ar[r]^{\epsilon_A} &    \Omega    \\
   A\times B \ar[u]^{1_A\times g} \ar[ur]_{f}                     }
  $$
      commutative.
\end{itemize}}
\end{definition}

Let us remain some facts from topos theory:
\begin{itemize}
  \item [(i)] A topos is a cartesian closed category.
  \item [(ii)] In a topos, all finite colimits exist.
  \item [(iii)] A topos has $(\texttt{Epi},\texttt{Mono})$ fs.
  \item [(iv)] Epimorphisms are stable under pullbacks.
  \item [(v)] If $\mathcal{C}$ is a topos, then for any object $A\in\mathcal{C}$ the comma category $\mathcal{C}\diagup A$ is also a topos.
  \item [(vi)] If $\mathbb{G}=(G,\epsilon,\delta)$ is a cartesian comonad (the functor $G$ preserves finite limits) in a topos $\mathcal{C}$, then the category of $\mathbb{G}$-coalgebras $\mathcal{C}_{\mathbb{G}}$ is a topos.
\end{itemize}

\begin{definition}[Definition V.1.1 \cite{MM}]\label{d13}
{\rm A \textit{Lawvere-Tierney topology} (\textit{LT-topology}) on a topos $\mathcal{C}$ is a morphism $k\colon\Omega\to\Omega$ for which the following diagrams commute:
$$
\xymatrix{\Omega \ar[dr]_{k} \ar[r]^{k} & \Omega \ar[d]^{k} & 1 \ar[l]_{\text{true}} \ar[dl]^{\text{true}}  \\
                & \Omega,             }
\
\
\xymatrix{\Omega\times\Omega \ar[d]_{k\times k} \ar[r]^{\wedge} & \Omega \ar[d]^{k}  \\
          \Omega\times\Omega                    \ar[r]_{\wedge} & \Omega,             }
$$
where $\wedge=\text{char}(\text{true}\times\text{true})$ (the morphism $\wedge$ is called the \textit{internal meet operation})}.
\end{definition}

The notion of LT-topology is closely related to the notion of a universal closure operator.

\begin{definition}[Lemma A.4.3.2 \cite{J}]\label{d14}
{\rm A map $c=c_A\colon\text{Sub}(A)\to\text{Sub}(A)$, $A\in\mathcal{C}$ is called a \textit{universal closure operator} if it satisfies the following axioms:
\begin{itemize}
  \item [(i)] $c(c(m))\simeq c(m)$;
  \item [(ii)] $m\le c(m)$;
  \item [(iii)] $c(m\cap n)\simeq c(m)\cap c(n)$;
  \item [(iv)] For any arrow $f\colon A\to B$, $f^*(c_B(m))\simeq c_A(f^*(m))$ holds, where $f^*\colon\mathcal{C}\diagup B\to\mathcal{C}\diagup A$ is the change-of-base functor.
\end{itemize}}
\end{definition}

LT-topologies are in bijective correspondence with universal closure operators (see Proposition V.1.1 \cite{MM}), i.e., $$c\mapsto k=\text{char}(c(\text{true}))\quad \text{and}\quad k\mapsto c(m)\simeq(k\text{char}(m))^*(\text{true}).$$ A monomorphism \(m\colon A\rightarrowtail B\) is called \textit{$k$-dense} (resp. \textit{$k$-closed}) if $1_B\simeq c(m)$ (resp. $m\simeq c(m)$).

From the definitions it follows immediately that $k$-dense (resp. $k$-closed) monomorphisms form a subcategory of $\texttt{Mono}$. The category of all $k$-dense (resp. $k$-closed) monomorphisms will be denoted by $\texttt{DnsMono}_k$ (resp. $\texttt{ClsMono}_k$). From Lemma A.4.3.3 \cite{J} it follows that $(\texttt{DnsMono}_k,\texttt{ClsMono}_k)$ is an fs in the category $\texttt{Mono}$, and moreover, dense and closed monomorphisms are stable under pullbacks. Now we can easily define the notion of a $k$-sheaf ($k$-separating object): an object $A$ is called a \textit{$k$-sheaf} (resp. \textit{$k$-separating object}) in the topos $\mathcal{C}$ if and only if $A$ is a $\texttt{DnsMono}_k$-local object (resp. $\texttt{DnsMono}_k$-separating object). The category of $k$-sheaves will be denoted by $\text{Sh}_k\ \mathcal{C}$, and the category of $k$-separating objects by $\text{Sep}_k\ \mathcal{C}$.

Finally, we note that the category of $k$-sheaves $\text{Sh}_k\ \mathcal{C}$ is itself a topos (Theorem V.2.5 \cite{MM}) and there is an adjunction
$\mathbf{a}\colon\mathcal{C}\leftrightarrows\text{Sh}_k\ \mathcal{C}\colon\iota$ such that the functor $\mathbf{a}$ preserves all finite limits (Theorem V.3.1 \cite{MM}). The functor $\mathbf{a}$ is called the \textit{sheafification functor}, and an adjunction between toposes whose left adjoint preserves all finite limits is called a \textit{geometric morphism}.

\subsection{Double factorization systems in a topos.}

Let us recall the notion of Bousfield localization (see Definition 7.20 \cite{JT}). More detailed information about Bousfield localizations can be found in \cite{Hr}. Let $(\mathcal{C}of_i,\mathcal{W}_i,\mathcal{F}ib_i)$, $i=1,2$ be two Quillen model structures in $\mathcal{C}$. The Quillen model structure $(\mathcal{C}of_2,\mathcal{W}_2,\mathcal{F}ib_2)$ is called a \textit{Bousfield localization of} $(\mathcal{C}of_1,\mathcal{W}_1,\mathcal{F}ib_1)$ if $\mathcal{C}of_1=\mathcal{C}of_2$ and $\mathcal{W}_1\subseteq\mathcal{W}_2$. We introduce an analogous definition for the case of dfs.

\begin{definition}\label{d15}
{\rm Let $(\mathcal{E}_i,\mathcal{J}_i,\mathcal{M}_i)$, $i=1,2$ be two dfs's in a category $\mathcal{C}$. We say that $(\mathcal{E}_2,\mathcal{J}_2,\mathcal{M}_2)$ is a \textit{Bousfield localization} of $(\mathcal{E}_1,\mathcal{J}_1,\mathcal{M}_1)$, if $\mathcal{J}_1\cdot\mathcal{E}_1=\mathcal{J}_2\cdot\mathcal{E}_2$ and $\mathcal{M}_1\cdot\mathcal{E}_1\subseteq\mathcal{M}_2\cdot\mathcal{E}_2$}.
\end{definition}

Let $k$ be an LT-topology on a topos $\mathcal{C}$, and let $\texttt{DnsMono}_k$ and $\texttt{ClsMono}_k$ be the classes of $k$-dense and $k$-closed monomorphisms respectively. From Lemma A4.3.3 \cite{J} and Proposition IV.6.2 \cite{MM}, it follows that $(\texttt{Epi},\texttt{DnsMono}_k,\texttt{ClsMono}_k)$ is a dfs in which the classes $\texttt{Epi}$, $\texttt{DnsMono}_k$ and $\texttt{ClsMono}_k$ are stable under pullbacks. In general, $(\texttt{Epi},\texttt{DnsMono}_k,\texttt{ClsMono}_k)$ does not define a qfs (see Remark \ref{r2}). The following definition is motivated by properties of the dfs $(\texttt{Epi},\texttt{DnsMono}_k,\texttt{ClsMono}_k)$.

\begin{definition}\label{d16}
{\rm A dfs $(\mathcal{E},\mathcal{J},\mathcal{M})$ on a topos $\mathcal{C}$ is called \textit{cartesian} if the following conditions hold:
\begin{itemize}
  \item [(i)] The classes $\mathcal{E}$ and $\mathcal{J}\cdot\mathcal{E}$ are stable under pullbacks;
  \item [(ii)] $\mathcal{M}\cdot\mathcal{J}\subseteq\texttt{Mono}$.
\end{itemize}}
\end{definition}

\begin{lemma}\label{l4}
For a dfs $(\mathcal{E},\mathcal{J},\mathcal{M})$ on a topos $\mathcal{C}$, the following conditions are equivalent:
\begin{itemize}
  \item [(i)] $(\mathcal{E},\mathcal{J},\mathcal{M})$ is a cartesian dfs;
  \item [(ii)] trivial fibrations and bifibrant morphisms are monomorphisms, and $(\mathcal{E},\mathcal{J},\mathcal{M})$ is stable under pullbacks.
\end{itemize}
\end{lemma}

\begin{proof}
(i)$\Longrightarrow$(ii) Property (ii) of Definition \ref{d16} is obvious. Since $(\mathcal{J}\cdot\mathcal{E},\mathcal{M})$ is an fs, the class $\mathcal{M}$ is stable under pullbacks. The stability of classes $\mathcal{E}$ and $\mathcal{J}\cdot\mathcal{E}$ follows from Definition \ref{d16}.

Consider the following commutative diagrams:
$$
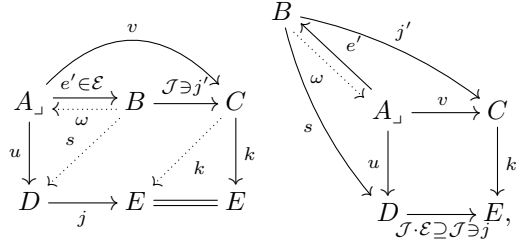
\xymatrix{ & &  \\
  A_{\lrcorner} \ar@/^2pc/[rr]^{v} \ar[d]_{u} \ar[r]<+0.5ex>^{e'\in\mathcal{E}}
  & B \ar@{.>}[l]<+0.5ex>^{\omega}\ar@{.>}[dl]_s \ar[r]^{\mathcal{J}\ni j'}  & C \ar@{.>}[dl]^k \ar[d]^k         \\
  D \ar[r]_{j} & E \ar@{=}[r]    &       E  }
\xymatrix{
   B \ar@/_/[ddr]_{s} \ar@/^/[drr]^{j'} \ar@{.>}[dr]<-0.5ex>_{\omega}                                    \\
   & A_{\lrcorner} \ar[lu]<-0.5ex>_{e'} \ar[d]_{u}     \ar[r]^{v}                             & C \ar[d]^{k}                     \\
   & D \ar[r]_{\mathcal{J}\cdot\mathcal{E}\supseteq\mathcal{J}\ni j}                                                 & E,                 }
$$
where the square (rectangle) is a pullback. This implies the equality $1_A=\omega e'$. Since $j'\in\texttt{Mono}$, then $\omega$ is a right-invertible monomorphism, i.e., an isomorphism. Consequently, $v\in\mathcal{J}$. The implication (ii)$\Longrightarrow$(i) is obvious.
\end{proof}

\begin{proposition}\label{p3}
A cartesian dfs $(\mathcal{E},\mathcal{J},\mathcal{M})$ on a topos $\mathcal{C}$ generates an \emph{LT}-topology $k$ and is a Bousfield localization of the dfs $(\emph{\texttt{Epi}},\emph{\texttt{DnsMono}}_k,\emph{\texttt{ClsMono}}_k)$.
\end{proposition}

\begin{proof}
Let $mje$ be a $(\mathcal{E},\mathcal{J},\mathcal{M})$-factorization of the subobject classifier $\text{true}\colon 1\rightarrowtail\Omega$ in $\mathcal{C}$. By Lemma \ref{l4}, $(\mathcal{E},\mathcal{J},\mathcal{M})$ is stable under pullbacks. In other words, for any morphism $f\colon A\to B$, the functor $f^*\colon\mathcal{C}\diagup B\to\mathcal{C}\diagup A$ preserves, up to isomorphism, the $(\mathcal{E},\mathcal{J},\mathcal{M})$-factorizations of morphisms in $\mathcal{C}$. In this case, the composition $(m\times m)(j\times j)(e\times e)$ is the $(\mathcal{E},\mathcal{J},\mathcal{M})$-factorization of the morphism $\text{true}\times\text{true}$. Let $k=\text{char}\ m$, then
\begin{align*}
(k\wedge)^*(\mathrm{true})&\simeq\wedge^*k^*(\mathrm{true})  \\
              &\simeq\wedge^*(m)  \\
              &\simeq m\times m   \\
              &\simeq k^*(\mathrm{true})\times k^*(\mathrm{true})  \\
              &\simeq (k\times k)^*(\mathrm{true\times true})  \\
              &\simeq (k\times k)^*(\wedge^*(\mathrm{true}))  \\
              &\simeq (\wedge(k\times k))^*(\mathrm{true}).\\
\end{align*}
Therefore, $k\wedge=\wedge(k\times k)$. The proof of $k=kk$ is similar. The equality $k\mathrm{true}=\mathrm{true}$ follows from the definition of \(k\). Thus, \(k\) is an LT-topology on \(\mathcal{C}\).

Let us consider the relationship between $(\mathcal{E},\mathcal{J},\mathcal{M})$ and $(\texttt{Epi},\texttt{DnsMono}_k,\texttt{ClsMono}_k)$. From the definitions of LT-topology and the class $\texttt{ClsMono}_k$, it follows that $\mathcal{M}=\texttt{ClsMono}_k$, and from this follows the equality $\mathcal{J}\cdot\mathcal{E}=\texttt{DnsMono}_k\cdot\texttt{Epi}$.
By Definition \ref{d16}, we have $\texttt{Epi}={}^\perp\texttt{Mono}\subseteq{}^\perp(\mathcal{M}\cdot\mathcal{J})=\mathcal{E}$. Thus, $\texttt{ClsMono}_k\cdot\texttt{Epi}\subseteq\mathcal{M}\cdot\mathcal{E}$. Therefore, $(\mathcal{E},\mathcal{J},\mathcal{M})$ is a Bousfield localization of the dfs $(\texttt{Epi},\texttt{DnsMono}_k,\texttt{ClsMono}_k)$.
\end{proof}

To conclude this section, let us consider the question of comparing LT-topologies. Recall that two LT-topologies $k_1$ and $k_2$ are in the relation $k_1\le k_2$ if and only if $k_1=\wedge(k_1,k_2)$. The relation $\le$ is a partial order on the set of all LT-topologies in a topos $\mathcal{C}$. It is known that the relation $\le$ can be characterized in terms of universal closure operators and sheaves. We can add a characterization of the relation $\le$ using cartesian dfs's. If $k_1$ and $k_2$ are two LT-topologies generated by cartesian dfs's $(\mathcal{C}of_1,\mathcal{W}_1,\mathcal{F}ib_1)$ and $(\mathcal{C}of_2,\mathcal{W}_2,\mathcal{F}ib_2)$ respectively, then
\begin{align*}
  k_1\le k_2 & \Longleftrightarrow \text{Sh}_{k_2}\mathcal{C}\subseteq\text{Sh}_{k_1}\mathcal{C}         \\
             & \Longleftrightarrow c_1(m)\le c_2(m),\quad m\in\texttt{Mono}                               \\
             & \Longleftrightarrow \texttt{DnsMono}_{k_1}\subseteq\texttt{DnsMono}_{k_2}                           \\
             & \Longleftrightarrow \mathcal{J}_1\cdot\mathcal{E}_1\subseteq\mathcal{J}_2\cdot\mathcal{E}_2                                      \\
\end{align*}

\begin{corollary}\label{c3}
Two cartesian dfs's $(\mathcal{E}_1,\mathcal{J}_1,\mathcal{M}_1)$ and $(\mathcal{E}_2,\mathcal{J}_2,\mathcal{M}_2)$ on a topos $\mathcal{C}$ define equal LT-topologies if and only if $\mathcal{J}_1\cdot\mathcal{E}_1=\mathcal{J}_2\cdot\mathcal{E}_2$ or $\mathcal{M}_1=\mathcal{M}_2$.
\end{corollary}

\subsection{Quillen adjunctions and the topos of coalgebras.}

We recall the notion (Definition in 5.7. \cite{DT}) of $(c,d)$-continuity (resp. $(c,d)$-preservation). Let us consider universal closure operators $c$ and $d$ of $\mathcal{C}$ and $\mathcal{D}$, respectively. A functor $F\colon\mathcal{C}\to\mathcal{D}$ is called \textit{$(c,d)$-continuous} (resp. \textit{$(c,d)$-preserving}) if $$F(c_B(f))\leq d_{FB}(F(f))\quad(\text{resp.}\quad F(c_B(f))\simeq d_{FB}(F(f))$$ holds for all $(\text{dom}(f),f)\in\mathcal{C}\diagup B$, $B\in\mathcal{C}$.

The following theorem establishes the connection between Quillen adjunctions and cartesian dfs's.

\begin{theorem}\label{th7}
Let $(\mathcal{E}_1,\mathcal{J}_1,\mathcal{M}_1)$ and $(\mathcal{E}_2,\mathcal{J}_2,\mathcal{M}_2)$ be cartesian dfs's on toposes $\mathcal{C}_1$ and $\mathcal{C}_2$ respectively, and let $F\colon\mathcal{C}_1\leftrightarrows\mathcal{C}_2\colon G$ be a Quillen adjunction. Then\emph{:}
\begin{itemize}
\item [(a)] if the left Quillen functor $F$ preserves finite limits \emph{(}the Quillen adjunction is a geometric morphism\emph{)}, then the right Quillen functor $G$ maps $k_2$-sheaves \emph{(}resp. $k_2$-separating objects\emph{)} to $k_1$-sheaves \emph{(}resp. $k_1$-separating objects\emph{);}
\item [(b)] the following conditions are equivalent\emph{:}
\begin{itemize}
  \item [(i)] $G(\mathcal{J}_2\cdot\mathcal{E}_2)\subseteq\mathcal{J}_1\cdot\mathcal{E}_1$\emph{;}
  \item [(ii)] $G$ is $(c_2,c_1)$-preserving\emph{;}
  \item [(iii)] $G$ is $(c_2,c_1)$-continuous\emph{;}
  \item [(iv)] $k_1\emph{\text{char}}(G(\emph{\text{true}}_{\mathcal{C}_2}))=\emph{\text{char}}(G(\emph{\text{true}}_{\mathcal{C}_2}))G(k_2)$,
\end{itemize}
\end{itemize}
where $k_1$ \emph{(}resp. $k_2$\emph{)} is the LT-topology on $\mathcal{C}_1$ \emph{(}resp. $\mathcal{C}_2$\emph{)} generated by the cartesian dfs $(\mathcal{E}_1,\mathcal{J}_1,\mathcal{M}_1)$ \emph{(}resp. $(\mathcal{E}_2,\mathcal{J}_2,\mathcal{M}_2)$\emph{)}, and $c_1$ \emph{(}resp. $c_2$\emph{)} is the universal closure operator corresponding to the LT-topology $k_1$ \emph{(}resp. $k_2$\emph{)}.
\end{theorem}

\begin{proof}
(a) By assumption, $F$ preserves finite limits, in particular, it maps monomorphisms to monomorphisms. Then
\begin{align*}
 \texttt{DnsMono}_{k_1}\subseteq\texttt{DnsMono}_{k_1}\cdot\texttt{Epi}=\mathcal{J}_1\cdot\mathcal{E}_1
          & \Longrightarrow F(\texttt{DnsMono}_{k_1})\subseteq\texttt{Mono}\cap F(\mathcal{J}_1\cdot\mathcal{E}_1)        \\
          & \Longrightarrow F(\texttt{DnsMono}_{k_1})\subseteq\texttt{Mono}\cap \mathcal{J}_2\cdot\mathcal{E}_2             \\
          & \Longrightarrow F(\texttt{DnsMono}_{k_1})\subseteq\texttt{Mono}\cap(\texttt{DnsMono}_{k_2}\cdot\texttt{Epi})      \\
          & \Longrightarrow F(\texttt{DnsMono}_{k_1})\subseteq\texttt{DnsMono}_{k_2}                                              \\
          & \Longrightarrow G(\text{Sh}_{k_2}\mathcal{C}_2)\subseteq\text{Sh}_{k_1}\mathcal{C}_1.                                       \\
\end{align*}
The last implication is essentially an equivalence and follows directly from the natural isomorphism of Hom-sets $\mathcal{C}_2(FA,B)\simeq\mathcal{C}_1(A,GB)$.

It is easy to verify that the diagonal monomorphism $\Delta_A\colon A\rightarrowtail A\times A$ is closed if and only if $A$ is a separating object (see the proof of implication (iv)$\Longrightarrow$(i) in Lemma A4.3.6 \cite{J} and the proof of implication (i)$\Longrightarrow$(ii) in Lemma V.3.3 \cite{MM}). Then we obtain
\begin{align*}
  \text{Sep}_{k_2}\mathcal{C}_2\ni A & \Longrightarrow\Delta_A\in\texttt{ClsMono}_{k_2}=\mathcal{M}_2 \\
                                     & \Longrightarrow\Delta_{GA}\simeq G(\Delta_A)\in\mathcal{M}_1=\texttt{ClsMono}_{k_1} \\
                                     & \Longrightarrow GA\in\text{Sep}_{k_1}\mathcal{C}_1.                                                 \\
\end{align*}

Let us prove the equivalence of the conditions listed in part (b).

(i)$\Longrightarrow$(ii) Consider an arbitrary morphism $f\colon A\to B$ in the topos $\mathcal{C}_2$ and let $mje$ be its $(\mathcal{E}_2,\mathcal{J}_2,\mathcal{M}_2)$-factorization. Since $\mathcal{M}_2\subseteq\mathcal{M}_2\cdot\mathcal{J}_2\subseteq\texttt{Mono}$, by definition of the universal closure operator, we have $c_2(f)\simeq m$. Without loss of generality, we assume that $f=c_2(f)je$. Let $c_1(G(f))j'e'$ be the $(\mathcal{E}_1,\mathcal{J}_1,\mathcal{M}_1)$-factorization of $Gf$. From the definition of the Quillen adjunction and condition (i) we have $G(je)\in G(\mathcal{J}_2\cdot\mathcal{E}_2)\subseteq\mathcal{J}_1\cdot\mathcal{E}_1$ and $G(c_2(f))\in G(\mathcal{M}_2)\subseteq\mathcal{M}_1$. From the uniqueness (up to isomorphism) of the $(\mathcal{E}_1,\mathcal{J}_1,\mathcal{M}_1)$-factorization of the morphism $Gf$ and the equality $G(c_2(f))G(j)G(e)=c_1(G(f))j'e'$ we obtain that $G(c_2(f))\simeq c_1(G(f))$.

(ii)$\Longrightarrow$(iii) Obviously.

(iii)$\Longrightarrow$(i) Taking into account the definition of universal closure operator, we obtain
\begin{align*}
  f\in\mathcal{J}_2\cdot\mathcal{E}_2 & \Longrightarrow c_2(f)\in\texttt{Iso} \\
                                      & \Longrightarrow\texttt{Iso}\ni G(c_2(f))\leq c_1(G(f)) \\
                                      & \Longrightarrow c_1(G(f))\in\texttt{Iso}                 \\
                                      & \Longrightarrow G(f)\in\mathcal{J}_1\cdot\mathcal{E}_1.
\end{align*}

(ii)$\Longleftrightarrow$(iv) Let $\tau=\text{char}(G(\text{true}_{\mathcal{C}_2}))$. The equivalence of conditions (ii) and (iv) follows directly from the relations
\begin{align*}
 (k_1\tau)^*(\text{true}_{\mathcal{C}_1}) & \simeq\tau^*k_1^*(\text{true}_{\mathcal{C}_1})  &  (\tau G(k_2))^*(\text{true}_{\mathcal{C}_1}) & \simeq  G(k_2)^*(\tau^*(\text{true}_{\mathcal{C}_1}))\\
                                        & \simeq\tau^*(c_1(\text{true}_{\mathcal{C}_1}))  && \simeq G(k_2)^*(G(\text{true}_{\mathcal{C}_2}))     \\
                                        & \simeq c_1(\tau^*(\text{true}_{\mathcal{C}_1})) && \simeq G(k_2^*(\text{true}_{\mathcal{C}_2}))         \\
                                        & \simeq c_1(G(\text{true}_{\mathcal{C}_2}))      && \simeq G(c_2(\text{true}_{\mathcal{C}_2})).            \\
\end{align*}
\end{proof}

\begin{remark}\label{r8}
{\rm Since LT-topologies are in one-to-one correspondence with universal closure operators, we can define $(k_1,k_2)$-continuous (resp. $(k_1,k_2)$-preserving) functor with respect to LT-topologies $k_1$ and $k_2$. More precisely, a functor $F\colon\mathcal{C}\to\mathcal{D}$ is \textit{$(k_1,k_2)$-continuous} (resp. \textit{$(k_1,k_2)$-preserving}), if $F$ is \textit{$(c_1,c_2)$-continuous} (resp. \textit{$(c_1,c_2)$-preserving}), where $c_j$ is universal closure operators corresponding to LT-topologies $k_j$, $j=1,2$. In the scope of Theorem \ref{th7}, we can conclude that the functor $F$ is $(k_1,k_2)$-continuous if and only if $F$ satisfies one of the equivalent conditions of item (b) of Theorem \ref{th7}.}
\end{remark}

Let $\mathcal{C}$ be a topos equipped with a cartesian dfs $(\mathcal{E},\mathcal{J},\mathcal{M})$ and a cartesian comonad $\mathbb{G}=(G,\epsilon,\delta)$, where the functor $G$ preserves all bifibrant morphisms and trivial fibrations, i.e., $G(\mathcal{J})\subseteq\mathcal{J}$ and $G(\mathcal{M})\subseteq\mathcal{M}$. By Proposition \ref{p2} and Proposition \ref{p3}, the cartesian dfs generates an LT-topology $k$ on $\mathcal{C}$, and the forgetful functor $G_{\mathbb{G}}$ induces a left-induced dfs $(\mathcal{E}_{\mathbb{G}},\mathcal{J}_{\mathbb{G}},\mathcal{M}_{\mathbb{G}})$ on the topos of $\mathbb{G}$-coalgebras $\mathcal{C}_{\mathbb{G}}$. Since $\mathcal{E}_{\mathbb{G}}=(G_{\mathbb{G}})^{-1}(\mathcal{E})$, $\mathcal{J}_{\mathbb{G}}=(G_{\mathbb{G}})^{-1}(\mathcal{J})$ and $\mathcal{M}_{\mathbb{G}}=(G_{\mathbb{G}})^{-1}(\mathcal{M})$, and the functor $G_{\mathbb{G}}$ reflects all finite limits, then $(\mathcal{E}_{\mathbb{G}},\mathcal{J}_{\mathbb{G}},\mathcal{M}_{\mathbb{G}})$ is a cartesian dfs on $\mathcal{C}_{\mathbb{G}}$. By Proposition \ref{p3}, the dfs $(\mathcal{E}_{\mathbb{G}},\mathcal{J}_{\mathbb{G}},\mathcal{M}_{\mathbb{G}})$ defines an LT-topology $k_{\mathbb{G}}$ on $\mathcal{C}_{\mathbb{G}}$.

Let $\tau=\text{char}(G(\text{true}))$ and let $m_{\mathbb{G}}$ be the equalizer in $\mathcal{C}$ of the arrows $1_{G\Omega}$, $G(\tau)\delta_{\Omega}$. Note that $1_{G\Omega}$ and $G(\tau)\delta_{\Omega}$ are arrows not only in $\mathcal{C}$ but also in $\mathcal{C}_{\mathbb{G}}$, so $m_{\mathbb{G}}$ is also their equalizer in $\mathcal{C}_{\mathbb{G}}$. The subobject classifier $\text{true}_{\mathbb{G}}$ in $\mathcal{C}_{\mathbb{G}}$ is defined as the unique arrow satisfying the equality $m_{\mathbb{G}}\text{true}_{\mathbb{G}}=G(\text{true})$ (see V.8 \cite{MM}).

A natural question arises: how are the LT-topologies $k$ and $k_{\mathbb{G}}$ related? To answer this question, we need the following lemma.

\begin{lemma}\label{l5}
Let $\mathbb{G}=(G,\epsilon,\delta)$ be a cartesian comonad on a topos $\mathcal{C}$ with an LT-topology $k$. If the functor $G$ is $(k,k)$-continuous, then the LT-topology $k$ can be extended to an LT-topology $\widetilde{k}$ on the topos of $\mathbb{G}$-coalgebras. Moreover, the forgetful functor $G_{\mathbb{G}}$ preserves and reflects the corresponding closed and dense monomorphisms.
\end{lemma}

\begin{proof}
There exists a unique arrow $\widetilde{k}\colon(\Omega_{\mathbb{G}},\omega)\to(\Omega_{\mathbb{G}},\omega)$ in the topos of $\mathbb{G}$-coalgebras making the diagram
$$
\xymatrix{
  \ar@{.>}[d]_{\widetilde{k}} \Omega_{\mathbb{G}} \ar[r]^{m_{\mathbb{G}}} & G\Omega \ar[d]_{Gk} \ar[r]^{\delta_{\Omega}} & G^2\Omega \ar[d]_{G^2k} \ar[r]^{G\tau} & G\Omega \ar[d]^{Gk} \\
  \Omega_{\mathbb{G}} \ar[r]_{m_{\mathbb{G}}} & G\Omega \ar[r]_{\delta_{\Omega}} & G^2\Omega  \ar[r]_{G\tau} & G\Omega }
$$
commute (the right square commutes by condition (iii) of item (b) of Theorem \ref{th7}).

To prove that $\widetilde{k}$ is an LT-topology, we need the explicit form of the internal meet operation in the topos of $\mathbb{G}$-coalgebras. We obtain:
\begin{align*}
  (\wedge(\tau\times\tau)\phi)^*(\text{true}) & \simeq \phi^*(\tau\times\tau)^*\wedge^*(\text{true})         \\
                                              & \simeq \phi^*(\tau\times\tau)^*(\text{true}\times\text{true}) \\
                                              & \simeq \phi^*(G(\text{true})\times G(\text{true}))             \\
                                              & \simeq G(\text{true}\times\text{true})                          \\
                                              & \simeq G(\wedge^*(\text{true}))                                  \\
                                              & \simeq G(\wedge)^*(G(\text{true}))                                \\
                                              & \simeq G(\wedge)^*(\tau^*(\text{true}))                            \\
                                              & \simeq (\tau G(\wedge))^*(\text{true})                              \\
                                              & \Longrightarrow \wedge(\tau\times\tau)\phi=\tau G(\wedge),           \\
\end{align*}
where $\phi\colon G(\Omega\times\Omega)\simeq G\Omega\times G\Omega$ is the natural isomorphism. This equality implies the existence of a unique arrow $\wedge_{\mathbb G}\colon(\Omega_{\mathbb G}\times\Omega_{\mathbb G},\phi^{-1}(\omega\times\omega))\to(\Omega_{\mathbb
G},\omega)$ in the category of $\mathbb{G}$-coalgebras making the diagram
$$\xymatrix{\Omega_{\mathbb{G}}\times\Omega_{\mathbb{G}} \ar@{.>}[ddd]_{\wedge_{\mathbb{G}}}\ar[rr]^{m_{\mathbb{G}}\times m_{\mathbb{G}}}
            && G\Omega\times G\Omega\ar[rr]^{\delta_{\Omega}\times\delta_{\Omega}} \ar@{-}[dd]_{\simeq}
            && G^2\Omega\times G^2\Omega \ar[rr]^{G(\tau)\times G(\tau)} \ar@{-}[d]_{\simeq}
            && G\Omega\times G\Omega\ar@{-}[d]_{\simeq} & \\
            &&  && G(G\Omega\times G\Omega) \ar@{-}[d]_{\simeq} \ar[rr]^{G(\tau\times\tau)}  &&  G(\Omega\times\Omega) \ar[d]^{G(\wedge)} \\
            && G(\Omega\times\Omega) \ar[d]_{G(\wedge)} \ar[rr]_{\delta_{\Omega\times\Omega}}
            && G^2(\Omega\times\Omega) \ar[r]_{G^2(\wedge)} & G^2(\Omega) \ar[r]_{G(\tau)}& \ar@{=}[d] G\Omega \\
            \Omega_{\mathbb{G}} \ar[rr]_{m_{\mathbb{G}}} && G(\Omega) \ar[rr]_{\delta_{\Omega}} & & G^2(\Omega) \ar[rr]_{G(\tau)} && G\Omega        }   $$
commute. We show that $\wedge_{\mathbb{G}}$ is the desired internal meet operation in the topos of $\mathbb{G}$-coalgebras. For this, it suffices to prove $\wedge_{\mathbb{G}}^*(\text{true}_{\mathbb{G}})\simeq\text{true}_{\mathbb{G}}\times\text{true}_{\mathbb{G}}$. Using the proved equality $\tau m_{\mathbb{G}}\wedge_{\mathbb{G}}=\tau G(\wedge)\phi(m_{\mathbb{G}}\times m_{\mathbb{G}})$, we obtain:
\begin{align*}
   \wedge_{\mathbb{G}}^*(\text{true}_{\mathbb{G}}) & \simeq \wedge_{\mathbb{G}}^*m_{\mathbb{G}}^*(G(\text{true}))                     \\
                                                   & \simeq \wedge_{\mathbb{G}}^*m_{\mathbb{G}}^*\tau^*(\text{true})                   \\
                                                   & \simeq (\tau m_{\mathbb{G}}\wedge_{\mathbb{G}})^*(\text{true})                     \\
                                                   & \simeq (\tau G(\wedge)\phi(m_{\mathbb{G}}\times m_{\mathbb{G}}))^*(\text{true})     \\
                                                   & \simeq (m_{\mathbb{G}}\times m_{\mathbb{G}})^*\phi^*G(\wedge)^*\tau^*(\text{true})   \\
                                                   & \simeq (m_{\mathbb{G}}\times m_{\mathbb{G}})^*\phi^*G(\wedge)^*(G(\text{true}))       \\
                                                   & \simeq (m_{\mathbb{G}}\times m_{\mathbb{G}})^*\phi^*(G(\wedge^*(\text{true})))         \\
                                                   & \simeq (m_{\mathbb{G}}\times m_{\mathbb{G}})^*\phi^*(G(\text{true}\times\text{true}))   \\
                                                   & \simeq (m_{\mathbb{G}}\times m_{\mathbb{G}})^*(G(\text{true})\times G(\text{true}))      \\
                                                   & \simeq \text{true}_{\mathbb{G}}\times\text{true}_{\mathbb{G}}.                            \\
\end{align*}

From the definitions of $\wedge_{\mathbb{G}}$ and $\widetilde{k}$ it follows directly that $\widetilde{k}$ is an LT-topology on the topos $\mathcal{C}_{\mathbb{G}}$.

Consider an arbitrary monomorphism \(m\colon(D,d)\rightarrowtail(A,s)\) in \(\mathcal{C}_{\mathbb{G}}\). Its characteristic arrow \(\text{char}_{\mathbb{G}}m\colon(A,s)\to(\Omega_{\mathbb G},\omega)\) in \(\mathcal{C}_{\mathbb{G}}\) satisfies \(m_{\mathbb{G}}\text{char}_{\mathbb{G}}m=G(\text{char}m)s\), with \((\text{char}_{\mathbb{G}}m)^*(\text{true}_{\mathbb{G}})\simeq m\) in \(\mathcal{C}\). Suppose \(m\) is a dense (resp. closed) monomorphism in \(\mathcal{C}\). In terms of characteristic arrows, density (resp. closedness) of \(m\) means \(k\text{char}m=\text{true}!_A\) (resp. \(k\text{char}m=\text{char}m\)). If \(m\) is dense (resp. closed) in \(\mathcal{C}\), then:
\begin{align*}
m_{\mathbb G}\widetilde{k}\text{char}_{\mathbb{G}}m&=G(k)m_{\mathbb G}\text{char}_{\mathbb{G}}m
&m_{\mathbb G}\widetilde{k}\text{char}_{\mathbb{G}}m&=G(k)m_{\mathbb G}\text{char}_{\mathbb{G}}m    \\
                                                    &=G(k)G(\text{char}\ m)s                     &&=G(k)G(\text{char}m)s                         \\
                                                    &=G(\mathrm{true})G(!_A)s                    &&=G(\text{char}m)s                               \\
                                                    &=m_{\mathbb G}\mathrm{true}_{\mathbb G} G(!_A)s &&=m_{\mathbb G}\text{char}_{\mathbb{G}}m      \\
                                                    & \Longrightarrow \widetilde{k}\text{char}_{\mathbb{G}}m=\mathrm{true}_{\mathbb G}G(!_A)s
                                                    &&\Longrightarrow \widetilde{k}\text{char}_{\mathbb{G}}m=\text{char}_{\mathbb{G}}m.               \\
\end{align*}
Thus, $m\colon(D,d)\rightarrowtail(A,s)$ is dense (resp. closed) in \(\mathcal{C}_{\mathbb{G}}\).

Conversely, suppose $m\colon(D,d)\rightarrowtail(A,s)$ is dense (resp. closed) in \(\mathcal{C}_{\mathbb{G}}\). As shown in V.8 \cite{MM}, \(\text{char}m=\tau G(\text{char}m)s\). Thus,
\begin{align*}
  k\text{char}m & = k\tau G(\text{char}m)s                                     & k\text{char}m & = k\tau G(\text{char}m)s \\
                & = \tau G(k)G(\text{char}m)s                                  && = \tau G(k)G(\text{char}m)s               \\
                & = \tau G(k)m_{\mathbb{G}}\text{char}_{\mathbb{G}}m           && = \tau G(k)m_{\mathbb{G}}\text{char}_{\mathbb{G}}m \\
                & = \tau m_{\mathbb G}\widetilde{k}\text{char}_{\mathbb{G}}m   && = \tau m_{\mathbb G}\widetilde{k}\text{char}_{\mathbb{G}}m \\
                & = \tau m_{\mathbb G}\mathrm{true}_{\mathbb G}G(!_A)s         && = \tau m_{\mathbb G}\text{char}_{\mathbb{G}}m               \\
                & = \tau G(\mathrm{true})G(!_A)s                               && = \tau G(\text{char}m)s                                      \\
                & = \mathrm{true}!_A                                           && =\text{char}m.                                               \\
\end{align*}
Therefore, $m\colon(D,d)\rightarrowtail(A,s)$ is dense (resp. closed) in $\mathcal{C}_{\mathbb{G}}$ if and only if $m$ is dense (resp. closed) in $\mathcal{C}$. In other words, we have proved that the forgetful functor $G_{\mathbb{G}}$ preserves and reflects dense and closed monomorphisms.
\end{proof}

The following proposition shows that under certain conditions, the LT-topology $k_{\mathbb{G}}$ is an extension $\widetilde{k}$ of the LT-topology $k$.

\begin{proposition}\label{p4}
Let $\mathcal{C}$ be a topos equipped with a cartesian dfs $(\mathcal{E},\mathcal{J},\mathcal{M})$ generating an LT-topology $k$, and a cartesian comonad $\mathbb{G}=(G,\epsilon,\delta)$ such that $G$ is $(k,k)$-continuous and preserves trivial fibrations and bifibrant morphisms. Then the LT-topology $k_{\mathbb{G}}$ generated by the left-induced cartesian dfs $(\mathcal{E}_{\mathbb{G}},\mathcal{J}_{\mathbb{G}},\mathcal{M}_{\mathbb{G}})$ on the topos of $\mathcal{C}_\mathbb{G}$ equals the extension $\widetilde{k}$ of the LT-topology $k$.
\end{proposition}

\begin{proof}
By Lemma \ref{l5}, there exists an LT-topology $\widetilde{k}$ on the topos $\mathcal{C}_{\mathbb{G}}$ extending the LT-topology $k$. Let $\texttt{Epi}_{\mathbb{G}}$ be the class of epimorphisms in $\mathcal{C}_{\mathbb{G}}$. The forgetful functor $G_{\mathbb{G}}$, being left adjoint to $R_{\mathbb{G}}$, preserves and reflects epimorphisms. From Proposition \ref{p2} and Lemma \ref{l5} it follows that $G_{\mathbb{G}}$ preserves and reflects cofibrations of cartesian dfs's, as well as dense and closed monomorphisms, i.e. $G^{-1}_{\mathbb{G}}(\mathcal{J}\cdot\mathcal{E})=\mathcal{J}_{\mathbb{G}}\cdot\mathcal{E}_{\mathbb{G}}$, $G^{-1}_{\mathbb{G}}(\texttt{DnsMono}_k)=\texttt{DnsMono}_{\widetilde{k}}$ and $G^{-1}_{\mathbb{G}}(\texttt{ClsMono}_k)=\texttt{ClsMono}_{\widetilde{k}}$. By Theorem 7 we have $G(\texttt{DnsMono}_k)\subseteq\texttt{DnsMono}_k$. Since $\texttt{Epi}\perp\texttt{DnsMono}_k$, by dual statement to Lemma \ref{l2} we obtain $\texttt{Epi}_{\mathbb{G}}\perp\texttt{DnsMono}_{\widetilde{k}}$ and $G^{-1}_{\mathbb{G}}(\texttt{DnsMono}_k\cdot\texttt{Epi})=\texttt{DnsMono}_{\widetilde{k}}\cdot\texttt{Epi}_{\mathbb{G}}$. For any morphism $e$ in $\mathcal{C}_{\mathbb{G}}$ we have:
\begin{align*}
  \mathcal{J}_{\mathbb{G}}\cdot\mathcal{E}_{\mathbb{G}}\ni e & \Longleftrightarrow e\in\mathcal{J}\cdot\mathcal{E}                         \\
                                                             & \Longleftrightarrow e\in\texttt{DnsMono}_k\cdot\texttt{Epi}                    \\
                                                             & \Longleftrightarrow e\in\texttt{DnsMono}_{\widetilde{k}}\cdot\texttt{Epi}_{\mathbb{G}}.
\end{align*}
By Corollary \ref{c3}, from $\mathcal{J}_{\mathbb{G}}\cdot\mathcal{E}_{\mathbb{G}}=\texttt{DnsMono}_{\widetilde{k}}\cdot\texttt{Epi}_{\mathbb{G}}$ we obtain $k_{\mathbb{G}}=\tilde{k}$.
\end{proof}

\section{Applications in equivariant topology}

\subsection{Double factorization systems in the categories $\mathbf{Unif}$, $\mathbf{Tych}$ and $\mathbf{Comp}$.}

In this chapter we will consider the following categories:
\begin{itemize}
  \item $\mathbf{Unif}$ is a category of uniform spaces and uniformly continuous maps between them;
  \item $\mathbf{Tych}$ is a category of Tychonoff spaces and continuous maps between them;
  \item $\mathbf{Comp}$ is a category of Hausdorff compact spaces and continuous maps between them;
  \item $\mathbf{CBUnif}$ is a category of totally bounded complete uniform spaces and uniformly continuous maps between them.
\end{itemize}

We recall the definitions of main classes of epimorphisms and monomorphisms. We restrict ourselves to the case of epimorphisms. By dual, everything said about certain classes of epimorphisms transfers to corresponding classes of monomorphisms. An epimorphism $e$ is called:
\begin{itemize}
  \item \textit{extremal} if from factorization $e=\mu\nu$, $\mu\in\texttt{Mono}$ it follows that $\mu\in\texttt{Iso}$;
  \item \textit{strong} if $e\in{}^\perp\texttt{Mono}$;
  \item \textit{regular} if $e$ is a coequalizer of some pair of arrows.
\end{itemize}

Classes of extremal, strong and regular epimorphisms will be denoted by $\texttt{ExEpi}$, $\texttt{StrEpi}$ and $\texttt{RegEpi}$ respectively. From Proposition 4.3.3 and Proposition 4.3.6 \cite{B} follows the relation $\texttt{RegEpi}\subseteq\texttt{StrEpi}\subseteq\texttt{ExEpi}$. If a category $\mathcal{C}$ has pullbacks, then classes of strong and extremal epimorphisms coincide, i.e. $\texttt{ExEpi}=\texttt{StrEpi}$ (see Lemma 1.3.2 \cite{J}). In particular, in a bicomplete category $\texttt{ExEpi}=\texttt{StrEpi}$ and $\texttt{ExMono}=\texttt{StrMono}$.

Let us consider the category $\textbf{Unif}$. General information about category $\textbf{Unif}$ can be found in \cite{I}. A surjective uniformly continuous map $h$ is called \textit{uniform quotient map} if for any map $f$, from the fact that $fh$ is uniformly continuous follows uniform continuity of $f$.

In \cite{Ku} was proved the existence of uniformly continuous map $h\colon(X,\mathcal{U})\to(X/\mathcal{U},\overline{\mathcal{U}})$, where $\mathcal{U}$ is some pseudouniformity on $X$, $X/\mathcal{U}$ is the quotient set consisting of equivalence classes $[x]=\bigcap_{\gamma\in\mathcal{U}}\text{St}(x,\gamma)$, $\overline{\mathcal{U}}$ is the uniformity on $X/\mathcal{U}$ whose base consists of covers of the form $\overline{\gamma}=\{(X/\mathcal{U})\setminus h(X\setminus V)\mid V\in\gamma\}$, $\gamma\in\mathcal{U}$ and $h$ is the canonical map assigning to each element $x$ its equivalence class $[x]$. The uniformity $\overline{\mathcal{U}}$ is called \textit{quotient uniformity of pseudouniformity} $\mathcal{U}$. The map $h$ satisfies the condition $\forall\gamma(\gamma\in\overline{\mathcal{U}}\Longleftrightarrow h^{-1}(\gamma)\in\mathcal{U})$.

Hereafter we will use subscripts $u$, $t$ and $c$ in notation of classes $\texttt{RegEpi}$, $\texttt{RegMono}$, $\texttt{Mono}$, $\texttt{Epi}$ etc. in categories $\mathbf{Unif}$, $\mathbf{Tych}$ and $\mathbf{Comp}$ respectively. In these categories, monomorphisms are precisely the injective mappings (for $\mathbf{Unif}$ this follows from Proposition 4 Ch.II \cite{I}; for $\mathbf{Tych}$ and $\mathbf{Comp}$ it follows from Example 1.7.7(b) Vol.1 \cite{B}).

\begin{lemma}\label{l6}
In the category $\mathbf{Unif}$ there exists a dfs $(\emph{\texttt{ExEpi}}_u,\emph{\texttt{Bim}}_u,\emph{\texttt{ExMono}}_u)$, where\emph{:}
\begin{itemize}
\item [(i)] $\emph{\texttt{Bim}}_u=$dense uniformly continuous injections\emph{;}
\item [(ii)] $\emph{\texttt{ExEpi}}_u=\emph{\texttt{RegEpi}}_u=\emph{\texttt{StrEpi}}_u=$uniform quotient mappings.
\end{itemize}
\end{lemma}

\begin{proof}
The category $\textbf{Unif}$ is bicomplete, locally small, and locally cosmall. According to Theorem 5.5 \cite{N} and its dual statement, the category $\textbf{Unif}$ admits factorization systems $(\texttt{ExEpi}_u,\texttt{Mono}_u)$ and $(\texttt{Epi}_u,\texttt{ExMono}_u)$.

Since $\texttt{Iso}_u\subseteq\texttt{ExEpi}_u\cap\texttt{Bim}_u\cap\texttt{ExMono}_u$ and the classes $\texttt{ExEpi}_u$, $\texttt{Bim}_u$, and $\texttt{ExMono}_u$ are closed under composition, the following relations hold $\texttt{Iso}_u\cdot\texttt{ExEpi}_u\subseteq\texttt{ExEpi}_u$, $\texttt{Iso}_u\cdot\texttt{Bim}_u\cdot\texttt{Iso}_u\subseteq\texttt{Bim}_u$, and $\texttt{ExMono}_u\cdot\texttt{Iso}_u\subseteq\texttt{ExMono}_u$.

Consider an arbitrary uniformly continuous mapping $f$. Let $\mu e$ be the $(\texttt{ExEpi}_u,\texttt{Mono}_u)$-factorization of $f$. Since $\mu$ has an $(\texttt{Epi}_u,\texttt{ExMono}_u)$-factorization $mj$, we obtain $f=mje$, where $e\in\texttt{ExEpi}_u$, $j\in\texttt{Bim}_u$, and $m\in\texttt{ExMono}_u$. Thus, every uniformly continuous mapping has an $(\texttt{ExEpi}_u,\texttt{Bim}_u,\texttt{ExMono}_u)$-factorization.

For any commutative rectangle
$$
\xymatrix{
  A \ar[d]_{u} \ar[r]^{\texttt{ExEpi}_u\ni e} & B \ar@{.>}[d]_{\omega}\ar@{.>}[dl]_s \ar[r]^{j\in\texttt{Bim}_u} & C \ar@{.>}[dl]^t \ar[d]^{v} \\
  D \ar[r]_{\texttt{Bim}_u\ni j'} & E \ar[r]_{m\in\texttt{ExMono}_u} & F   }
$$
there exists a unique arrow $\omega$ that splits it into two commutative squares. In turn, each square splits into commutative triangles via arrows $s$ and $t$. The uniqueness of $s$ and $t$ is obvious. We finally conclude that $(\texttt{ExEpi}_u,\texttt{Bim}_u,\texttt{ExMono}_u)$ is a factorization system in $\mathbf{Unif}$.

(i) From Proposition 13 Ch.I \cite{I}, it follows that every epimorphism $f\colon(X,\mathcal{U}_X)\to(Y,\mathcal{U}_Y)$ in $\mathbf{Unif}$ has a dense image, i.e., $\overline{f(X)}=Y$. Thus, every bimorphism in $\mathbf{Unif}$ is a dense uniformly continuous injection. The converse follows from the fact that any two uniformly continuous mappings that coincide on a dense subset are equal.

(ii) We have previously noted the inclusion $\texttt{RegEpi}_u\subseteq\texttt{StrEpi}_u\subseteq\texttt{ExEpi}_u$. Consider an arbitrary extremal epimorphism $e\colon(X,\mathcal{U}_X)\to(Y,\mathcal{U}_Y)$ in $\mathbf{Unif}$. The kernel pair of $e$ consists of the restrictions $\pi_i\colon X\times_Y X\to X$, $i=1,2$ of the projections $\text{pr}_i\colon X\times X\to X$, $i=1,2$ to the uniform space $(X\times_Y X,\mathcal{U}_{X\times_Y X})$, where $X\times_Y X=\{(x,y)\mid e(x)=e(y)\}$ and $\mathcal{U}_{X\times_Y X}=(\mathcal{U}_X\times\mathcal{U}_X)\cap(X\times_Y X)$ is the trace of the uniformity $\mathcal{U}_X\times\mathcal{U}_X$. Since every extremal epimorphism in $\mathbf{Unif}$ is a surjective uniformly continuous mapping, we obtain that $\texttt{coeq}(\pi_1,\pi_2)=\nu$ is a surjective uniformly continuous mapping. There exists a unique uniformly continuous mapping $e'$ satisfying $e=e'\nu$. Note that $e(x)=e(y)\Longleftrightarrow\nu(x)=\nu(y)$. We obtain $e'\in\texttt{Mono}_u\Longrightarrow e'\in\texttt{Iso}_u\Longrightarrow e\in\texttt{RegEpi}_u$. Consequently, $\texttt{RegEpi}_u=\texttt{StrEpi}_u=\texttt{ExEpi}_u$.

Let $e$ be a uniform quotient mapping and $e=me'$, where $m\in\texttt{Mono}_u$. Since $e$ is surjective, then $m$ is bijective. Therefore, $m^{-1}e=e'$ is uniformly continuous. From the definition of uniform quotient mapping, it follows that $m^{-1}$ is uniformly continuous, which means $m$ is an isomorphism in $\mathbf{Unif}$. Hence, $e\in\texttt{ExEpi}_u$. The converse follows from the fact that every regular epimorphism is a uniform quotient mapping.
\end{proof}

The category $\mathbf{Unif}$ is closely related to the category $\mathbf{Tych}$. The connection between these categories is realized via functors $T$ and $F$. Let us recall their definitions.

\textit{Topologization functor} $T\colon\mathbf{Unif}\to\mathbf{Tych}$ assigns to each uniform space $(X,\mathcal{U})$ a Tychonoff space $(X,\tau(\mathcal{U}))$, with topology $\tau(\mathcal{U})$ induced by uniformity $\mathcal{U}$. To each uniformly continuous map $f\colon(X,\mathcal{U})\to(Y,\mathcal{V})$, the functor $T$ assigns the same map $f\colon(X,\tau(\mathcal{U}))\to(Y,\tau(\mathcal{V}))$, considered as continuous in topologies $\tau(\mathcal{U})$ and $\tau(\mathcal{V})$ induced by uniformities $\mathcal{U}$ and $\mathcal{V}$ respectively.

\textit{Finest uniformity functor} $F\colon\mathbf{Tych}\to\mathbf{Unif}$. The subbase of the finest uniformity $\mathcal{U}^m_X$ consists of all open covers of space $(X,\tau)$. The functor $F$ assigns to Tychonoff space $(X,\tau)$ a uniform space $(X,\mathcal{U}^m_X)$, and to continuous map $f\colon(X,\tau_X)\to(Y,\tau_Y)$ assigns a uniformly continuous map $f\colon(X,\mathcal{U}^m_X)\to(Y,\mathcal{U}^m_Y)$ (Proposition 21 Ch.I \cite{I}).

\begin{lemma}\label{l7}
In the category $\mathbf{Tych}$, there exists a dfs $(\emph{\texttt{ExEpi}}_t,\emph{\texttt{Bim}}_t,\emph{\texttt{ExMono}}_t)$, where\emph{:}
\begin{itemize}
  \item [(i)] $\emph{\texttt{Bim}}_t=$ dense continuous injections\emph{;}
  \item [(ii)] $\emph{\texttt{RegEpi}}_t=\emph{\texttt{ExEpi}}_t=\emph{\texttt{StrEpi}}_t$.
\end{itemize}
Moreover, the adjunction $F\colon\mathbf{Tych}\leftrightarrows\mathbf{Unif}\colon T$ is a Quillen adjunction.
\end{lemma}

\begin{proof}
By Proposition 19 Ch.I \cite{I}, each cover $\gamma\in\mathcal{U}$ has an open in topology $\tau(\mathcal{U})$ star-refinement. Therefore, the identity map $\epsilon_{(X,\mathcal{U})}=\text{id}\colon FT(X,\mathcal{U})\to(X,\mathcal{U})$ is uniformly continuous. From the definitions of $F$ and $T$, it follows that $TF = 1_{\mathbf{Tych}}$. Consequently, there exists an adjunction $F\colon\mathbf{Tych}\leftrightarrows\mathbf{Unif}\colon T$ with $\eta = 1_{1_{\mathbf{Tych}}}\colon 1_{\mathbf{Tych}} \Longrightarrow TF$ as the unit and $\epsilon\colon FT \Longrightarrow 1_{\mathbf{Unif}}$ as the counit (in this case, $F$ is called \textit{left adjoint and right inverse to} $T$).

From the statement dual to Exercise 4 \S4 Ch.4 \cite{M}, it follows that there exists a complete coreflective subcategory $\mathcal{D}$ in $\mathbf{Unif}$ and an isomorphism $H\colon\mathbf{Tych}\simeq\mathcal{D}$ such that $F=KH$, where $K\colon\mathcal{D}\hookrightarrow\mathbf{Unif}$ is the inclusion. From the statements dual to Proposition 3.5.3 and Proposition 3.5.4 Vol.1 \cite{B}, it follows that the category $\mathbf{Tych}$ is bicomplete. In particular, in the category $\mathbf{Tych}$, strong monomorphisms (resp. strong epimorphisms) coincide with external monomorphisms (resp. external epimorphisms), that is, $\texttt{StrMono}_t=\texttt{ExMono}_t$ (resp. $\texttt{StrEpi}_t=\texttt{ExEpi}_t$).

By Proposition 4.3.9 Vol.1 \cite{B}, we obtain $T(\texttt{ExMono}_u)\subseteq\texttt{ExMono}_t$, $F(\texttt{ExEpi}_t)\subseteq\texttt{ExEpi}_u$, and $F(\texttt{Epi}_t)\subseteq\texttt{Epi}_u$. Therefore, $F(\texttt{Bim}_t)\subseteq\texttt{Bim}_u$, and hence $F\colon\mathbf{Tych}\leftrightarrows\mathbf{Unif}\colon T$ is a Quillen adjunction.

(i) Every continuous dense injection is a bimorphism in $\mathbf{Tych}$. On the other hand, a uniformly continuous dense injection is a continuous dense injection in the topologies induced by uniformities. Therefore, $T(\texttt{Bim}_u)\subseteq\texttt{Bim}_t$. Taking into account the relation $F(\texttt{Bim}_t)\subseteq\texttt{Bim}_u$, we finally obtain $\texttt{Bim}_t=T(\texttt{Bim}_u)$=dense continuous injections.

(ii) Let $e\in\texttt{ExEpi}_t$. Then $Fe\in\texttt{ExEpi}_u=\texttt{RegEpi}_u$, that is, there exists a pair of uniformly continuous mappings $f,g$ such that $\texttt{coeq}(f,g)=Fe$. Since $\epsilon_F=1_F$ and $F(\texttt{coeq}(Tf,Tg))\simeq\texttt{coeq}(FTf,FTg)$, there exist arrows $\zeta$ and $\xi$ satisfying the equalities $e=\zeta\texttt{coeq}(Tf,Tg)$ and $F(\texttt{coeq}(Tf,Tg))=\xi Fe$. Therefore, $e\simeq\texttt{coeq}(Tf,Tg)\in\texttt{RegEpi}_t$. Finally, we obtain $\texttt{ExEpi}_t=\texttt{RegEpi}_t=\texttt{StrEpi}_t$.
\end{proof}

\begin{lemma}\label{l8}
The restriction of the dfs \emph{$(\texttt{ExEpi}_t,\texttt{Bim}_t,\texttt{ExMono}_t)$} from the category $\mathbf{Tych}$ on the subcategory $\mathbf{Comp}$ induces the dfs \emph{$(\texttt{RegEpi}_c,\texttt{Iso}_c,\texttt{RegMono}_c)$}. The following diagram commutes\emph{:}
$$
\xymatrix{\mathbf{Tych} \ar@<-0.5ex>[d]_{\beta} \ar@<-0.5ex>[r]_F & \mathbf{Unif} \ar@<-0.5ex>[l]_T \ar@<-0.5ex>[d]_{\beta^u} \\
\mathbf{Comp} \ar@{_{(}->}@<-0.5ex>[u]_{\iota} \ar@<-0.5ex>[r]_{F^c} & \mathbf{CBUnif} \ar@<-0.5ex>[l]_{T^c} \ar@{_{(}->}@<-0.5ex>[u]_{\iota^u}, }
$$
where\emph{:}
\begin{itemize}
\item $F^c\colon\mathbf{Comp}\leftrightarrows\mathbf{CBUnif}\colon T^c$ is the restriction of $F\colon\mathbf{Tych}\leftrightarrows\mathbf{Unif}\colon T$ and induces an isomorphism between $\mathbf{Comp}$ and $\mathbf{CBUnif}$\emph{;}
\item $\beta^u\colon\mathbf{Unif}\to\mathbf{CBUnif}$ is the Samuel compactification\emph{;}
\item $\beta\colon\mathbf{Tych}\to\mathbf{Comp}$ is the Stone-Cech compactification\emph{;}
\item $\iota$ and $\iota^u$ are the corresponding embedding functors.
\end{itemize}
All adjunctions in the diagram are Quillen adjunctions.
\end{lemma}

\begin{proof}
Since $\beta\colon\mathbf{Tych}\to\mathbf{Comp}$ is a reflector, it follows from Proposition 3.5.3 and Proposition 3.5.4 Vol.1 \cite{B} that the category $\mathbf{Comp}$ is bicomplete, and hence locally small and locally cosmall.

From Proposition 4.3.9 Vol.1 \cite{B}, it follows that the embedding functor $\iota$ preserves external and regular monomorphisms, i.e., $\texttt{ExMono}_c\subseteq\texttt{ExMono}_t$ and $\texttt{RegMono}_c\subseteq\texttt{RegMono}_t$. Any continuous map $f$ in $\mathbf{Comp}$ can be represented as a composition $me$, where $e$ is a surjection and $m$ is a closed embedding. Hence, $\texttt{ExMono}_t\cap\mathbf{Comp}\subseteq\{\text{closed continuous embeddings}\}\cap\mathbf{Comp}$. Let $Z$ be a closed subset of a compact space $Y$. Consider the space $Y/Z$ obtained by collapsing $Z$ to a point. Let $h\colon Y\to Y/Z$ be the quotient map, and let $g\colon Y\to Y/Z$ be the constant map mapping all points of $Y$ to the equivalence class $Z\in Y/Z$. Clearly, the embedding $m\colon Z\hookrightarrow Y$ satisfies $hm=gm$. If $f\colon X\to Y$ is a continuous map in $\mathbf{Comp}$ such that $hf=gf$, then the image $f(X)$ is a (closed) subset of $Z$. Consequently, there exists a unique continuous map $k$ such that $f=mk$. It follows that $\texttt{RegMono}_c=\{\text{closed continuous embeddings}\}\cap\mathbf{Comp}$, and therefore $\texttt{ExMono}_t\cap\mathbf{Comp}=\texttt{RegMono}_c$.

From Example 4.3.10.f Vol.1 \cite{B}, it follows that $\texttt{Epi}_c=\texttt{RegEpi}_c=\texttt{StrEpi}_c=\texttt{ExEpi}_c$. Consequently, $\texttt{ExEpi}_t\cap\mathbf{Comp}\subseteq\texttt{RegEpi}_c$. Let $\texttt{ExEpi}_c\ni f=\mu\nu$, where $\mu\in\texttt{Mono}_t$. Decompose $\nu$ as $je$, where $e$ is a surjection and $j$ is an embedding. Since the continuous image of a compact space is compact, we have $\mu j\in\texttt{Bim}_c=\texttt{Mono}_c\cap\texttt{Epi}_c=\texttt{Mono}_c\cap\texttt{StrEpi}_c=\texttt{Iso}_c\subseteq\texttt{Iso}_t\Longrightarrow\mu\in\texttt{Iso}_t$. Hence, $\texttt{ExEpi}_c\subseteq\texttt{ExEpi}_t\cap\mathbf{Comp}$. We conclude that $\texttt{ExEpi}_t\cap\mathbf{Comp}=\texttt{RegEpi}_c$.

By Lemma \ref{l7}, the class $\texttt{Bim}_t$ comprises precisely the dense continuous injections. Since a dense continuous injection from a compact space into a compact space is a homeomorphism, then $\texttt{Bim}_t\cap\mathbf{Comp}=\texttt{Iso}_c$.

The above results imply that the $(\texttt{ExEpi}_t,\texttt{Bim}_t,\texttt{ExMono}_t)$-factorization of an arbitrary morphism $f$ in $\mathbf{Comp}$ coincides with its $(\texttt{RegEpi}_c,\texttt{Iso}_c,\texttt{RegMono}_c)$-factorization in $\mathbf{Comp}$, and $\beta\colon\mathbf{Tych}\leftrightarrows\mathbf{Comp}\colon\iota$ is a Quillen adjunction.

Since a uniform space is compact if and only if it is both totally bounded and complete, the restriction of the adjunction $F\colon\mathbf{Tych}\leftrightarrows\mathbf{Unif}\colon T$ to $\mathbf{Comp}$ and $\mathbf{CBUnif}$ yields an adjunction $F^c\colon\mathbf{Comp}\leftrightarrows\mathbf{CBUnif}\colon T^c$, which is an isomorphism of categories and hence a Quillen adjunction (the fs in $\mathbf{CBUnif}$ is defined as the image of $(\texttt{RegEpi}_c,\texttt{Iso}_c,\texttt{RegMono}_c)$ under the isomorphism $F^c$).

By Theorem 32 Ch.II \cite{I}, there exists an adjunction $\beta^u\colon\mathbf{Unif}\leftrightarrows\mathbf{CBUnif}\colon\iota^u$, where $\beta^u$ is the Samuel compactification functor. The equality $T\iota^u=\iota T^c$ implies that (up to natural isomorphism) $\beta^uF=F^c\beta$. The embedding functor $\iota^u$, being right adjoint, maps regular monomorphisms to regular monomorphisms. In the category $\mathbf{CBUnif}$, bimorphisms are precisely the isomorphisms. Thus, $\iota^u$ maps bimorphisms to bimorphisms. Consequently, $\beta^u\colon\mathbf{Unif}\leftrightarrows\mathbf{CBUnif}\colon\iota^u$ is a Quillen adjunction.
\end{proof}

\subsection{Double factorization systems in the categories $\mathbf{Tych}^{\mathbb{H}^t}$, $\mathbf{Unif}^{\mathbb{H}^u}$ and $\mathbf{Comp}^{\mathbb{H}^c}$.}

Throughout this section, unless stated otherwise, we consider an arbitrary fixed compact topological group $G$. Let us consider the category of all (small) topological spaces $\mathbf{Top}$ and continuous maps between them. We define a monad in $\mathbf{Top}$ according to Section 1.1 \S1 \cite{V}. Consider the functor $H\colon\mathbf{Top}\to\mathbf{Top}$ as the cartesian product of $G$, that is, $X \overset{H}{\mapsto} G\times X$, $f \overset{H}{\mapsto} 1_G\times f$. The unit $\eta$ and multiplication $\mu$ of the monad are defined as follows: $$ X\ni x \overset{\eta_X}{\mapsto} (e,x)\in G\times X\quad \text{and} \quad G\times(G\times X)\ni (g_1,(g_2,x)) \overset{\mu_X}{\mapsto} (g_1g_2,x)\in G\times X.$$
The commutativity of diagrams from Definition \ref{d10} and continuity of maps $\eta_X$ and $\mu_X$ are obvious. Therefore, $\mathbb{H}=(H,\eta,\mu)$ is a monad in $\mathbf{Top}$.

Now let $\mathcal{B}$ be one of the categories $\mathbf{Tych}$, $\mathbf{Unif}$ or $\mathbf{Comp}$. To define a monad in a similar way, two conditions must be satisfied (see Section 4.1 \S4 \cite{V}):
\begin{itemize}
  \item [(M1)] for each object $X\in\mathcal{B}$ the product $G\times X$ is an object of $\mathcal{B}$;
  \item [(M2)] for each object $X\in\mathcal{B}$ the maps $\eta_X$ and $\mu_X$ are morphisms in $\mathcal{B}$.
\end{itemize}

Since $G$ is a compact group, the restriction $H^t$ of the functor $H$ to the category $\mathbf{Tych}$ satisfies condition (M1), and the restrictions $\eta^t$ of the unit $\eta$ and $\mu^t$ of the multiplication $\mu$ to objects of $\mathbf{Tych}$ satisfy condition (M2). Therefore, $\mathbb{H}^t=(H^t,\eta^t,\mu^t)$ is a monad in $\mathbf{Tych}$. Moreover, compactness of $G$ is a necessary and sufficient condition for the analogous construction of the monad $\mathbb{H}^c=(H^c,\eta^c,\mu^c)$ in $\mathbf{Comp}$.

We define the monad $\mathbb{H}^u$ in $\mathbf{Unif}$ as follows. The functor $H^u\colon\mathbf{Unif}\to\mathbf{Unif}$ is defined by the rule: $H^u(X,\mathcal{U})=(G,\mathcal{U}_G)\times(X,\mathcal{U})$ and $H^u(f)=1_G\times f$, where $\mathcal{U}_G$ is the unique totally bounded and complete uniformity on the group $G$ that generates its original (group) topology. The maps $\eta^u$ and $\mu^u$ are defined in a similar way to $\eta$ and $\mu$. The functor $H^u$ and the maps $\eta^u$, $\mu^u$ satisfy conditions (M1) and (M2) (the multiplication in a compact group is a uniformly continuous map). Therefore, $\mathbb{H}^u=(H^u,\eta^u,\mu^u)$ is a monad in $\mathbf{Unif}$.

The total boundedness and completeness of the uniformity $\mathcal{U}_G$ on $G$ are necessary and sufficient conditions for the analogous definition of the monad $\mathbb{H}^b=(H^b,\eta^b,\mu^b)$ in $\mathbf{CBUnif}$.

Let $s$ be one of the indices $t$, $u$ or $c$. Then $\mathcal{B}^{\mathbb{H}^s}$ denotes the category of $\mathbb{H}^s$-algebras in $\mathcal{B}$.

\begin{proposition}\label{p5}
The category $\mathcal{B}^{\mathbb{H}^s}$ is \emph{(}small\emph{)} cocomplete.
\end{proposition}

\begin{proof}
\textit{Case $s=u$}.
Let $\{((X_j,\mathcal{U}_j),\alpha_j)\}_{j\in J}$ be a set of objects in $\mathbf{Unif}^{\mathbb{H}^u}$. The action $\coprod_{j\in J}\alpha_j$ of the group $G$ on the disjoint union of sets $\coprod_{j\in J}X_j$ is uniquely determined by the commutative diagram
$$
\xymatrix{G\times X_j \ar[d]_{\alpha_j} \ar@{^{(}->}[r]^{1_G\times\sigma_j} & G\times\coprod\limits_{j\in J}X_j \ar[d]^{\coprod\limits_{j\in J}\alpha_j}                     \\
             X_j       \ar@{^{(}->}[r]_{\sigma_j}                           & \coprod\limits_{j\in J}X_j,             }
$$
where $\sigma_j\colon X_j\hookrightarrow\coprod_{j\in J}X_j$, $j\in J$ are the canonical embeddings. We define the uniformity $\coprod_{j\in J}\mathcal{U}_j$ on $\coprod_{j\in J}X_j$ as follows: $$\coprod_{j\in J}\mathcal{U}_j\ni\gamma\Longleftrightarrow\forall j(j\in J\Longrightarrow\sigma^{-1}_j(\gamma)\in\mathcal{U}_j).$$ The fact that $\coprod_{j\in J}\mathcal{U}_j$ is a uniformity on $\coprod_{j\in J}X_j$ follows from Proposition 2 Ch.II \cite{I}. In this case, the action $\coprod_{j\in J}\alpha_j\colon(G,\mathcal{U}^m_G)\times(\coprod_{j\in J}X_j,\coprod_{j\in J}\mathcal{U}_j)\to(\coprod_{j\in J}X_j,\coprod_{j\in J}\mathcal{U}_j)$ is uniformly continuous, and the set $\{\sigma_j\colon((X_j,\mathcal{U}_j),\alpha_j)\hookrightarrow((\coprod_{j\in J}X_j,\coprod_{j\in J}\mathcal{U}_j),\coprod_{j\in J}\alpha_j)$ is the coproduct in $\mathbf{Unif}^{\mathbb{H}^u}$.

The idea of the proof of existence of coequalizers in $\mathbf{Unif}^{\mathbb{H}^u}$ is close to the case of $\mathbf{Tych}^G$ (see Theorem 1 \cite{ME}). Consider two morphisms $f_1,f_2\colon((X,\mathcal{U}_X),\alpha)\to((Y,\mathcal{U}_Y),\beta)$ in $\mathbf{Unif}^{\mathbb{H}^u}$. Let $\mathfrak{S}$ be the set of pseudouniformities $\mathcal{U}$ on $Y$ satisfying the following conditions:
\begin{itemize}
  \item [(a)] $\mathcal{U}\subseteq\mathcal{U}_Y$;
  \item [(b)] the restriction $\beta|_{(G,\mathcal{U}^m_G)\times(Y,\mathcal{U})}$ of the action $\beta$ is a uniformly continuous map from $(G,\mathcal{U}^m_G)\times(Y,\mathcal{U})$ to $(Y,\mathcal{U})$;
  \item [(c)] for any $x\in X$ and any $\gamma\in\mathcal{U}$, there exists $V\in\gamma$ such that $\{f_1(x),f_2(x)\}\subseteq V$.
\end{itemize}

The set $\mathfrak{S}$ is non-empty, as it contains the pseudouniformity to which the cover $\{Y\}$ belongs. Consider the supremum $\bigvee\mathfrak{S}=\bigvee_{\mathcal{U}\in\mathfrak{S}}\mathcal{U}$ of pseudouniformities in $\mathfrak{S}$. A base for $\bigvee\mathfrak{S}$ consists of covers of the form $\bigwedge_{j=1}^m\gamma_j$, $\gamma_j\in\mathcal{U}_j\in\mathfrak{S}$. Condition (a) for $\bigvee\mathfrak{S}$ is obvious. For any $\gamma_j\in\mathcal{U}_j\in\mathfrak{S}$ there exist $\gamma'_j\in\mathcal{U}_G$ and $\gamma''_j\in\mathcal{U}_j$ such that $(\bigwedge_{j=1}^m\gamma'_j)\times(\bigwedge_{j=1}^m\gamma''_j)=\bigwedge_{j=1}^m(\gamma'_j\times\gamma''_j)\succ\bigwedge_{j=1}^m\beta^{-1}(\gamma_j)$. Consequently, $\beta|_{(G,\mathcal{U}^m_G)\times(Y,\bigvee\mathfrak{S})}$ is a uniformly continuous map from $(G,\mathcal{U}^m_G)\times(Y,\bigvee\mathfrak{S})$ to $(Y,\bigvee\mathfrak{S})$. For any $x\in X$ and any $\bigwedge_{j=1}^m\gamma_j$, $\gamma_j\in\mathcal{U}_j\in\mathfrak{S}$ there exist $V_j\in\gamma_j\in\mathcal{U}_j$ such that $\{f_1(x),f_2(x)\}\subseteq\bigcap_{j=1}^mV_j\in\bigwedge_{j=1}^m\gamma_j\in\bigvee\mathfrak{S}$. Thus, $\bigvee\mathfrak{S}$ is the largest element in $\mathfrak{S}$.

For any $x\in X$ and arbitrary $\gamma\in\bigvee\mathfrak{S}$ there exists $V\in\gamma$ such that $\{f_1(x),f_2(x)\}\subseteq V$, and hence $f_1(x)\in\text{St}(f_2(x),\gamma)$. Since the choice of $\gamma$ is arbitrary, we obtain $f_1(x)\in[f_2(x)]=\bigcap_{\gamma\in\bigvee\mathfrak{S}}\text{St}(f_2(x),\gamma)$. From the definition of the uniform quotient map $h\colon(Y,\bigvee\mathfrak{S})\to(Y/\bigvee\mathfrak{S},\overline{\bigvee\mathfrak{S}})$ follows the equality $hf_1=hf_2$.

Consider the map $\eta^g_Y\colon Y\to G\times Y$ defined by $\eta^g_Y(y)=(g,y)$. Obviously, the maps $\eta^g_Y\colon(Y,\bigvee\mathfrak{S})\to(G,\mathcal{U}^m_G)\times(Y,\bigvee\mathfrak{S})$, $g\in G$ are uniformly continuous. Since $\beta^g=\beta\eta^g_Y$, then $\beta^g\colon(Y,\bigvee\mathfrak{S})\to(Y,\bigvee\mathfrak{S})$, $g\in G$ are uniformly continuous maps. Therefore, there exist uniformly continuous maps $\overline{\beta}^g\colon(Y/\bigvee\mathfrak{S},\overline{\bigvee\mathfrak{S}})\to(Y/\bigvee\mathfrak{S},\overline{\bigvee\mathfrak{S}})$ satisfying $\overline{\beta}^gh=h\beta^g$, $g\in G$.

The map $\overline{\beta}(g,[y])=\overline{\beta}^g(h(y))$ is a uniformly continuous action of the group $(G,\mathcal{U}_G)$ on the quotient space $(Y/\bigvee\mathfrak{S},\overline{\bigvee\mathfrak{S}})$. The fact that $\overline{\beta}$ satisfies the action axioms follows from the equalities $\overline{\beta}^gh=h\beta^g$, $g\in G$. Let us prove uniform continuity of $\overline{\beta}$. For any $\overline{\gamma}\in\overline{\bigvee\mathfrak{S}}$, we have $(1_G\times h)^{-1}(\overline{\beta}^{-1}(\overline{\gamma}))=\beta^{-1}(h^{-1}(\overline{\gamma}))\in\mathcal{U}_G\times\bigvee\mathfrak{S}$. Therefore, there exist covers $\gamma_1\in\mathcal{U}_G$ and $\gamma_2\in\bigvee\mathfrak{S}$ such that $\gamma_1\times\gamma_2\succ\beta^{-1}(h^{-1}(\overline{\gamma}))$. From the relation
\begin{align*}
\mathcal{U}_G\times\overline{\bigvee\mathfrak{S}}\ni\gamma_1\times\overline{\gamma}_2   & \succ\gamma_1\times h(\gamma_2)           \\
                             & =(1_G\times h)(\gamma_1\times\gamma_2)                                                                \\
                             & \succ(1_G\times h)(\beta^{-1}(h^{-1}(\overline{\gamma})))                                              \\
                             & =(1_G\times h)((1_G\times h)^{-1}(\overline{\beta}^{-1}(\overline{\gamma})))                            \\
                             & =\overline{\beta}^{-1}(\overline{\gamma})                                                                \\
\end{align*}
follows uniform continuity of $\overline{\beta}$. By definition of $\overline{\beta}$ we obtain that $h\colon((Y,\mathcal{U}_Y),\beta)\to((Y/\bigvee\mathfrak{S},\overline{\bigvee\mathfrak{S}}),\overline{\beta})$ is a uniformly continuous equivariant map.

Let $f\colon((Y,\mathcal{U}_Y),\beta)\to((Z,\mathcal{U}_Z),\delta)$ be an equivariant uniformly continuous map satisfying $ff_1=ff_2$. The map $f$ defines a pseudouniformity $f^{-1}(\mathcal{U}_Z)$ on $Y$. A base for $f^{-1}(\mathcal{U}_Z)$ consists of covers of the form $\bigwedge_{j=1}^mf^{-1}(\gamma_j)$, $\gamma_j\in\mathcal{U}_Z$. Since $f$ is uniformly continuous, $f^{-1}(\mathcal{U}_Z)$ satisfies condition (a). For any $\gamma\in\mathcal{U}_Z$ there exist covers $\gamma_1\in\mathcal{U}_G$ and $\gamma_2\in\mathcal{U}_Z$ such that $\gamma_1\times\gamma_2\succ\delta^{-1}(\gamma)$. From the relation
\begin{align*}
\mathcal{U}_G\times f^{-1}(\mathcal{U}_Z)\ni\gamma_1\times f^{-1}(\gamma_2) & = (1_G\times f)^{-1}(\gamma_1\times\gamma_2) \\
                                                                              & \succ (1_G\times f)^{-1}(\delta^{-1}(\gamma)) \\
                                                                              & =\beta^{-1}(f^{-1}(\gamma))                     \\
\end{align*}
follows that $f^{-1}(\mathcal{U}_Z)$ satisfies condition (b). For any $x\in X$ and arbitrary cover $\gamma\in\mathcal{U}_Z$ there exists $V\in\gamma$ such that $f(f_1(x))=f(f_2(x))\in V$. Thus, $\{f_1(x),f_2(x)\}\subseteq f^{-1}(V)\in f^{-1}(\mathcal{U}_Z)$. Therefore, $f^{-1}(\mathcal{U}_Z)$ satisfies condition (c), hence $f^{-1}(\mathcal{U}_Z)\in\mathfrak{S}$. From item (5) \cite{Ku} follows the existence and uniqueness of an equivariant uniformly continuous map $f'\colon((Y/\bigvee\mathfrak{S},\overline{\bigvee\mathfrak{S}}),\overline{\beta})\to((Z,\mathcal{U}_Z),\delta)$ satisfying $f=f'\circ h$. Consequently, we obtain $\texttt{coeq}(f_1,f_2)=h$ in $\mathbf{Unif}^{\mathbb{H}^u}$.

\textit{Case $s=t$}. The existence of coproducts in the category $\mathbf{Tych}^{\mathbb{H}^t}$ follows directly from the existence of coproducts in $\mathbf{Unif}^{\mathbb{H}^u}$. It suffices to carry out the same reasoning with the necessary modifications.

Consider two morphisms $f_1,f_2\colon(X,\alpha)\to(Y,\beta)$ in the category $\mathbf{Tych}^{\mathbb{H}^t}$. The coequalizer $\texttt{coeq}(f_1,f_2)=e\colon Y\twoheadrightarrow Z$ in $\mathbf{Tych}$ is a surjective map ($e\in\texttt{ExEpi}_t\Longrightarrow e$ is a surjection). Moreover, $e$ satisfies the following condition: for any map $f$, from the fact that $fe$ is continuous follows continuity of $f$. Then for any $g\in G$ there exists a unique continuous map $\delta^g$ such that $\delta^ge=e\beta^g$. Thus, the action $\delta\colon G\times Z\to Z$ of $G$ on $Z$ is defined. Since $G$ is compact, the evaluation map $\text{ev}\colon G\times Z^G\to Z$ is continuous ($Z^G$ is endowed with the compact-open topology). There exist a unique map $\ulcorner\delta\urcorner\colon Y\to Z^G$ satisfying $\text{ev}(1_G\times\ulcorner\delta\urcorner)=\delta$. We obtain that $e\beta$ is continuous $\Longrightarrow\ulcorner\delta\urcorner e$ is continuous $\Longrightarrow\ulcorner\delta\urcorner$ is continuous $\Longrightarrow\delta$ is continuous. The fact that $\texttt{coeq}(f_1,f_2)=e\colon(Y,\beta)\twoheadrightarrow(Z,\delta)$ in $\mathbf{Tych}^{\mathbb{H}^t}$ is verified directly.

\textit{Case $s=c$}. In \cite{V} it is proved that the category $\mathbf{Comp}^G$ is cocomplete. The difference between the category $\mathbf{Comp}^{\mathbb{H}^c}$ and the category $\mathbf{Comp}^G$ is in the acting group: in $\mathbf{Comp}^{\mathbb{H}^c}$ the acting group $G$ is compact, while in $\mathbf{Comp}^G$ it is arbitrary. Therefore, the category $\mathbf{Comp}^{\mathbb{H}^c}$ is cocomplete.
\end{proof}

For $\mathbb{H}^s$-algebras we will use the notation $\mathbb{X}^s$. There exists an adjunction $F^{\mathbb{H}^s}\colon\mathcal{B}\leftrightarrows\mathcal{B}^{\mathbb{H}^s}\colon V^{\mathbb{H}^s}$, where the right adjoint $V^{\mathbb{H}^s}$ is the forgetful functor and the left adjoint $F^{\mathbb{H}^s}$ is the free algebra functor.

\begin{theorem}\label{th8}
Let $\mathbb{H}^s=(H^s,\eta^s,\mu^s)$ be the monad in $\mathcal{B}$ described above. The following conditions hold\emph{:}
\begin{itemize}
  \item [(i)] $H^s$ preserves the dfs \emph{($\texttt{RegEpi}_s,\texttt{Bim}_s,\texttt{RegMono}_s$);}
  \item [(ii)] the category $\mathcal{B}^{\mathbb{H}^s}$ is \emph{(}small\emph{)} bicomplete\emph{;}
  \item [(iii)] $V^{\mathbb{H}^s}$ generates a right-induced dfs $(\emph{\texttt{RegEpi}}^{\mathbb{H}^s},\emph{\texttt{Bim}}^{\mathbb{H}^s},\emph{\texttt{RegMono}}^{\mathbb{H}^s})$\emph{;}
  \item [(iv)] $V^{\mathbb{H}^s}$ preserves and reflects the corresponding $\emph{\texttt{Epi}}$-local objects, i.e. for any $\mathbb{X}^s\in\mathcal{B}^{\mathbb{H}^s}$ the equivalence holds\emph{:} $$\mathbb{X}^s\in\mathcal{B}^{\mathbb{H}^s}_{\emph{\texttt{Epi}}^{\mathbb{H}^s}}\Longleftrightarrow V^{\mathbb{H}^s}(\mathbb{X}^s)\in\mathcal{B}_{\emph{\texttt{Epi}}_s};$$
  \item [(v)] The following diagram is commutative
$$
\xymatrix{
  \mathbf{Tych}^{\mathbb{H}^t} \ar@<-0.5ex>[dd]_{V^{\mathbb{H}^t}} \ar@<-0.5ex>[rd]_{\beta^{tc}} \ar@<-0.5ex>[rr]_{F^{tu}}
  && \mathbf{Unif}^{\mathbb{H}^u} \ar@<-0.5ex>[dd]_{V^{\mathbb{H}^u}} \ar@<-0.5ex>[ll]_{T^{tu}} \ar@<-0.5ex>[ld]_{L^{uc}}   \\
  & \mathbf{Comp}^{\mathbb{H}^c} \ar@{_{(}->}@<-0.5ex>[ul]_{\iota^{tc}} \ar@<-0.5ex>[ur]_{R^{uc}} \ar@<-0.5ex>[dd]_(0.299999){V^{\mathbb{H}^c}}
             \\
  \mathbf{Tych}  \ar@<-0.5ex>[uu]_{F^{\mathbb{H}^t}} \ar@<-0.5ex>[rd]_{\beta} \ar@<-0.5ex>[rr]_(0.62){F}
  && \mathbf{Unif} \ar@<-0.5ex>[ll]_(0.32){T} \ar@<-0.5ex>[ld]_{L^u} \ar@<-0.5ex>[uu]_{F^{\mathbb{H}^u}}                                \\
  & \mathbf{Comp}   \ar@<-0.5ex>[uu]_(0.699999){F^{\mathbb{H}^c}} \ar@<-0.5ex>[ru]_{R^u} \ar@{_{(}->}@<-0.5ex>[ul]_{\iota},          }
$$
$$\emph{\text{Fig. 1}}$$
where\emph{:}
\end{itemize}
\begin{itemize}
  \item $F^{tu}\colon\mathbf{Tych}^{\mathbb{H}^t}\leftrightarrows\mathbf{Unif}^{\mathbb{H}^u}\colon T^{tu}$ is the extension of $F\colon\mathbf{Tych}\leftrightarrows\mathbf{Unif}\colon T$\emph{;}
  \item $\beta^{tc}\colon\mathbf{Tych}^{\mathbb{H}^t}\leftrightarrows\mathbf{Comp}^{\mathbb{H}^c}\colon\iota^{tc}$ is the extension of $\beta\colon\mathbf{Tych}\leftrightarrows\mathbf{Comp}\colon\iota$\emph{;}
  \item $L^{uc}\colon\mathbf{Unif}^{\mathbb{H}^u}\leftrightarrows\mathbf{Comp}^{\mathbb{H}^c}\colon R^{uc}$ is the extension of $L^u\colon\mathbf{Unif}\leftrightarrows\mathbf{Comp}\colon R^u$, where $R^u=i^uF^c$ and $L^u=T^c\beta^u$.
\end{itemize}
All adjunctions in the commutative diagram are Quillen adjunctions.
\end{theorem}

\begin{proof}
(i) Since the category $\mathcal{B}$ is complete, we have the equality $\texttt{StrMono}_s=\texttt{ExMono}_s$. From the definition of $H^s$ and Lemmas \ref{l6}, \ref{l7} and \ref{l8} follow the relations $H^s(\texttt{ExMono}_s)\subseteq\texttt{ExMono}_s$ and $H^s(\texttt{Bim}_s)\subseteq\texttt{Bim}_s$. To prove $H^s(\texttt{ExEpi}_s)\subseteq\texttt{ExEpi}_s$, we consider each case separately.

\textit{Case $s=u$}. Lemma \ref{l6} implies that $\texttt{ExEpi}_u$ coincides with uniform quotient maps. By Theorem 1 \cite{Hu} the direct product of two uniform quotient maps is a uniform quotient map. In particular, if $f$ is a uniform quotient map, then the direct product $1_G \times f$ is also a uniform quotient map. Consequently, we obtain $H^u(\texttt{ExEpi}_u)\subseteq\texttt{ExEpi}_u$.

\textit{Case $s=t$}.
By assumption, $G$ is compact. Then there exists an adjunction $H^t=G\times(-)\colon\mathbf{Tych}\leftrightarrows\mathbf{Tych}\colon(-)^G$,
where $X^G$ is the space of continuous functions from $G$ to $X$ endowed with the compact-open topology (the fact that $X^G$ is a Tychonoff space follows from Theorem 3.4.15 \cite{E}). From Proposition 4.3.9 Vol.1 \cite{B} and Lemma \ref{l7} it follows that $H^t(\texttt{ExEpi}_t)\subseteq\texttt{ExEpi}_t$.

\textit{Case $s=c$}. By Lemma \ref{l8} the dfs $(\texttt{RegEpi}_c,\texttt{Iso}_c,\texttt{RegMono}_c)$ is a restriction of the dfs $(\texttt{ExEpi}_t,\texttt{Bim}_t,\texttt{ExMono}_t)$. The restriction of the functor $H^t$ to the category $\mathbf{Comp}$  yields the functor $H^c$. Consequently, $H^c$ preserves the dfs $(\texttt{RegEpi}_c,\texttt{Iso}_c,\texttt{RegMono}_c)$.

(ii) The completeness of the category $\mathcal{B}^{\mathbb{H}^s}$ follows from Proposition 4.3.1 Vol.2 \cite{B}, and the cocompleteness follows from Proposition \ref{p5}.

(iii) and (iv) follow from Proposition \ref{p1}.

(v) If in the diagram in Fig.~1 all faces commute with respect to the right adjoint functors, then they will also commute (up to natural isomorphism) with respect to the left adjoint functors. Therefore, it is sufficient to consider commutativity only with respect to the right adjoints. The commutativity of the upper triangle follows from the commutativity of the lower one and the side faces.

From the equalities $V^{\mathbb{H}^t}T^{tu}=TV^{\mathbb{H}^u}$, $V^{\mathbb{H}^t}\iota^{tc}=\iota V^{\mathbb{H}^c}$ and Theorem 4.5.6. Vol.2 \cite{B} we obtain the existence of adjunctions $F^{tu}\colon\mathbf{Tych}^{\mathbb{H}^t}\leftrightarrows\mathbf{Unif}^{\mathbb{H}^u}\colon T^{tu}$, $\beta^{tc}\colon\mathbf{Tych}^{\mathbb{H}^t}\leftrightarrows\mathbf{Comp}^{\mathbb{H}^c}\colon\iota^{tc}$ and the commutativity of the left and back faces of the diagram in Fig. 1. Consider the following commutative diagram
$$
\xymatrix{
\mathbf{Unif}^{\mathbb{H}^u}     \ar@<-0.5ex>[d]_{V^{\mathbb{H}^u}}      \ar@<-0.5ex>[r]_{\beta^{ub}}
& \mathbf{CBUnif}^{\mathbb{H}^b} \ar@<-0.5ex>[d]_{V^{\mathbb{H}^b}}      \ar@{_{(}->}@<-0.5ex>[l]_{\iota^{ub}}   \ar@<-0.5ex>[r]_{T^{cb}}
& \mathbf{Comp}^{\mathbb{H}^c}   \ar@<-0.5ex>[d]_{V^{\mathbb{H}^c}}      \ar@<-0.5ex>[l]_{F^{cb}}                      \\
\mathbf{Unif}                    \ar@<-0.5ex>[u]_{F^{\mathbb{H}^u}}      \ar@<-0.5ex>[r]_{\beta^u}
& \mathbf{CBUnif}                \ar@<-0.5ex>[r]_{T^c}                   \ar@{_{(}->}@<-0.5ex>[l]_{\iota^u}                   \ar@<-0.5ex>[u]_{F^{\mathbb{H}^b}}
& \mathbf{Comp}                  \ar@<-0.5ex>[u]_{F^{\mathbb{H}^c}}      \ar@<-0.5ex>[l]_{F^c}.   }
$$
The adjunction $F^{cb}\colon\mathbf{Comp}^{\mathbb{H}^c}\leftrightarrows\mathbf{CBUnif}^{\mathbb{H}^b}\colon T^{cb}$, as an extension of the isomorphism of categories $F^c\colon\mathbf{Comp}\leftrightarrows\mathbf{CBUnif}\colon T^c$, itself defines an isomorphism of categories ($F^{cb}$ and $T^{cb}$ are mutually inverse functors). In particular, $\mathbf{CBUnif}^{\mathbb{H}^b}$ is cocomplete and the right square in the diagram commutes. From Theorem 4.5.6. Vol.2 \cite{B} follows the existence of a left adjoint $\beta^{ub}\colon\mathbf{Unif}^{\mathbb{H}^u}\to\mathbf{CBUnif}^{\mathbb{H}^b}$ to the embedding functor $\iota^{ub}$ and commutativity of the left square. Therefore, the compositions $R^u=\iota^uF^c$ and $L^u=T^c\beta^u$ define an adjunction $L^u\colon\mathbf{Unif}\leftrightarrows\mathbf{Comp}\colon R^u$, and the compositions $R^{uc}=i^{ub}F^{cb}$ and $L^{uc}=T^{cb}\beta^{ub}$ define an adjunction $L^{uc}\colon\mathbf{Unif}^{\mathbb{H}^u}\leftrightarrows\mathbf{Comp}^{\mathbb{H}^c}\colon R^{uc}$.
In particular, the right face in the diagram in Fig. 1 commutes. Directly from definitions follows commutativity of the lower triangle, i.e. $\iota=\iota T^cF^c=T\iota^uF^c=TR^u$. To complete the proof, it suffices to refer to Lemmas \ref{l6}, \ref{l7}, \ref{l8} and Corollary \ref{c2}.
\end{proof}

\

\textsc{National Research Nuclear University "MEPhI", 31 Kashirskoe Highway, Moscow, 115409, Russia}

\textit{Email address}: binom00@yandex.ru

\end{document}